\def\DateTime{09/April/1998, 11:15AM}
\def\Version{Version 2.0}

\documentclass[leqno,12pt]{amsart}
\usepackage{amssymb}
\usepackage{amscd}
\usepackage{latexsym}
\usepackage{verbatim}

\theoremstyle{plain}
\newtheorem{Theorem}{Theorem}[section]
\newtheorem{Proposition}[Theorem]{Proposition}
\newtheorem{Lemma}[Theorem]{Lemma}
\newtheorem{Corollary}[Theorem]{Corollary}
\newtheorem{Claim}{Claim}[Theorem]
\newtheorem{Sublemma}[Claim]{Sublemma}

\theoremstyle{definition}

\newtheorem{Remark}[Theorem]{Remark}

\renewcommand{\theTheorem}{\arabic{section}.\arabic{Theorem}}
\renewcommand{\theClaim}{\arabic{section}.\arabic{Theorem}.\arabic{Claim}}
\renewcommand{\theequation}{\arabic{section}.\arabic{Theorem}.\arabic{Claim}}

\def\rom{\textup}
\newcommand{\ZZ}{{\mathbb{Z}}}
\newcommand{\QQ}{{\mathbb{Q}}}
\newcommand{\RR}{{\mathbb{R}}}
\newcommand{\CC}{{\mathbb{C}}}

\newcommand{\OO}{{\mathcal{O}}}
\newcommand{\codim}{\operatorname{codim}}

\newcommand{\adeg}{\widehat{\operatorname{deg}}}
\newcommand{\Image}{\operatorname{Im}}
\newcommand{\Supp}{\operatorname{Supp}}
\newcommand{\Spec}{\operatorname{Spec}}
\newcommand{\Sym}{\operatorname{Sym}}
\newcommand{\Proj}{\operatorname{Proj}}
\newcommand{\Sing}{\operatorname{Sing}}
\newcommand{\Ker}{\operatorname{Ker}}
\newcommand{\Rat}{\operatorname{Rat}}
\newcommand{\Br}{\mathbf{Br}}
\newcommand{\Cycle}{Z}
\newcommand{\aCycle}{\widehat{Z}}
\newcommand{\aBCycle}{\widehat{Z}_{B}}
\newcommand{\aBrCycle}{\widehat{Z}_{\Br}}
\newcommand{\aDCycle}{\widehat{Z}_D}
\newcommand{\Chow}{\operatorname{CH}}
\newcommand{\aChow}{\widehat{\operatorname{CH}}}
\newcommand{\aBChow}{\widehat{\operatorname{CH}}_{B}}
\newcommand{\aBrChow}{\widehat{\operatorname{CH}}_{\Br}}
\newcommand{\aDChow}{\widehat{\operatorname{CH}}_D}
\newcommand{\Div}{\operatorname{Div}}
\newcommand{\aDiv}{\widehat{\operatorname{Div}}}
\newcommand{\aBDiv}{\widehat{\operatorname{Div}}_{B}}
\newcommand{\aBrDiv}{\widehat{\operatorname{Div}}_{\Br}}

\newcommand{\Pic}{\operatorname{Pic}}
\newcommand{\aPic}{\widehat{\operatorname{Pic}}}
\newcommand{\aBPic}{\widehat{\operatorname{Pic}}_{B}}
\newcommand{\aBrPic}{\widehat{\operatorname{Pic}}_{\Br}}
\newcommand{\aDPic}{\widehat{\operatorname{Pic}}_D}
\newcommand{\poscone}{\widehat{C}_{+\kern-.09em+}}
\newcommand{\wposcone}{\widehat{C}_{+}}

\newcommand{\chernch}{\operatorname{ch}}
\newcommand{\acherncl}{\widehat{{c}}}
\newcommand{\achernch}{\widehat{\operatorname{ch}}}
\newcommand{\alev}{\operatorname{a.e.}}
\newcommand{\zero}{\operatorname{div}}
\newcommand{\rank}{\operatorname{rk}}
\newcommand{\Proof}{{\sl Proof.}\quad}
\newcommand{\QED}{{\unskip\nobreak\hfil\penalty50\quad\null\nobreak\hfil
{$\Box$}\parfillskip0pt\finalhyphendemerits0\par\medskip}}
\newcommand{\rest}[2]{\left.{#1}\right\vert_{{#2}}}

\begin{document}

\title[Intersection pairing for arithmetic cycles]%
{Intersection pairing for arithmetic cycles \\
with degenerate Green currents}
\author{Atsushi Moriwaki}
\address{Department of Mathematics, Faculty of Science,
Kyoto University, Kyoto, 606-01, Japan}
\email[Atsushi Moriwaki]{moriwaki@kusm.kyoto-u.ac.jp}
\date{\DateTime, (\Version)}

\maketitle


\section*{Introduction}
\renewcommand{\theTheorem}{\Alph{Theorem}}

As we indicated in our paper \cite{KMRB},
the standard arithmetic Chow groups introduced by
Gillet-Soul\'{e} \cite{GSArInt} are rather restricted to
consider arithmetic analogues of geometric problems.
In this note, we would like to propose a suitable extension
of the arithmetic Chow group of codimension one, in which
the Hodge index theorem still holds 
as in papers \cite{Fa}, \cite{Hr} and \cite{MoHI}.

Let $X \to \Spec(\ZZ)$ be a regular arithmetic variety with
$d = \dim X_{\QQ}$.
As we defined in \cite{KMRB},
$\aDChow^p(X)$ is a group, consisting
of pairs $(Z, g)$ with cycles $Z$ of codimension $p$ on $X$  and
currents $g$ of type $(p-1,p-1)$ on $X(\CC)$,
modulo arithmetical rational equivalence.
It seems to be impossible to give
a natural ring structure on the graded module
$\aDChow^{*}(X)_{\QQ}$.
In \cite[\S2.3]{KMRB}, we showed that 
$\aDChow^{*}(X)_{\QQ}$ has, however,
a natural scalar product 
of the arithmetic Chow ring $\aChow^{*}(X)_{\QQ}$,
namely, a module structure of $\aDChow^{*}(X)_{\QQ}$ 
over $\aChow^{*}(X)_{\QQ}$
as a generalization of \cite[Theorem~4.2.3]{GSArInt}.
In this note, we will introduce
suitable subgroups $\aChow_B^p(X)$ of $\aDChow^p(X)$ and
$\aChow_{B'}^q(X)$ of $\aDChow^q(X)$
such that $\aChow_B^p(X)$ and $\aChow_{B'}^q(X)$ have a natural paring
\[
\aChow_B^p(X) \otimes \aChow_{B'}^q(X) \to \aDChow^{p+q}(X)_{\QQ}.
\]
In the following, we would like to explain how to construct
$\aBChow^1(X)$, for example.

\medskip
We denote by $C^{\infty}(X(\CC), \RR)$ (resp. $L^2_{1,loc}(X(\CC), \RR)$)
the set of all real valued $C^{\infty}$-functions
(resp. locally square integrable functions with all weak partial derivatives in $L^2_{loc}$) 
on $X(\CC)$.
A key point to get $\aBChow^1(X)$ is to fix an abelian group $B$ with
$C^{\infty}(X(\CC), \RR) \subseteq B \subseteq L^{2}_{1,loc}(X(\CC), \RR)$.
This abelian group $B$ is 
called a class of degeneration of Green functions.
Fixing the class $B$ of degeneration,
an arithmetic $B$-divisor on $X$ is defined to be 
a pair $(D, g)$ such that
$D$ is a divisor on $X$, and that
there are a Green function $f$ for $D(\CC)$ and $\phi \in B$
with $g = f + \phi$.
We denote by $\aBCycle^1(X)$ the set of all
arithmetic $B$-divisors on $X$, and define
$\aBChow^1(X)$ to be $\aBCycle^1(X)$ modulo
arithmetic linear equivalence.
Using the Dirichlet form of $L^2_{1,loc}(X(\CC), \RR)$
together with the usual star product,
if $(D_1, g_1), (D_2, g_2) \in \aBCycle^1(X)$,
and $D_1$ and $D_2$ have no common component, then we can
define the star product $g_1 * g_2$ as well as the usual one
(cf. \S\ref{subsec:green:fun:via:class}). In this way,
we have a homomorphism
\[
\aBChow^1(X) \otimes \aBChow^1(X) \to \aDChow^2(X).
\]
as desired.

\medskip
Assuming $X$ is projective over $\ZZ$,
let $(H, k)$ be an arithmetically ample Hermitian line bundle on $X$,
namely,
(1) $H$ is ample,
(2) the Chern form $c_1(H, k)$ is 
positive definite on the infinite fiber $X(\CC)$, and
(3) there is a positive integer $m_0$ such that,
for any integer $m \geq m_0$, $H^0(X, H^m)$
is generated by the set
$\left\{ s \in H^0(X, H^m) \mid 
\Vert s \Vert_{\sup} < 1 \right\}$.
Then, the following is one of main results of this note,
which is a generalization of \cite{Fa}, \cite{Hr} and \cite{MoHI}.

\begin{Theorem}[cf. Corollary~\ref{cor:standard:conj:div}]
Let us consider a homomorphism
\[
L : \aChow^p_D(X)_{\QQ} \to \aChow^{p+1}_D(X)_{\QQ}
\]
given by $L(x) = \acherncl_1(\overline{H}) \cdot x$. Then, we have the following.
\begin{enumerate}
\renewcommand{\labelenumi}{(\arabic{enumi})}
\item
$L^{d-1} : \aChow^1_D(X)_{\QQ} \to \aChow^{d}_D(X)_{\QQ}$
is injective.

\item
If $x \in \aChow^1_{B}(X)_{\QQ}$, $x \not= 0$, and
$L^d(x) = 0$, then $\adeg(L^{d-1}(x) \cdot x) < 0$.
\end{enumerate}
\end{Theorem}

On an arithmetic surface, Bost \cite{Bo} also constructed the same
intersection pairing independently
for an arithmetic analogue of the Lefschetz theorem.
Our motivation is, however, different.
In this note, we introduce the special class $\Br$,
which is called the birational class of degeneration,
arising from birational geometry. 
Namely, a locally integrable function $\phi$ on $X(\CC)$
belongs to the class $\Br$ if and only if
there are a proper birational morphism $\mu : Y \to X(\CC)$
of smooth algebraic schemes over $\CC$, a divisor $D'$ on $Y$, and
a Green function $g$ for $D'$ such that
$\mu_*(D') = 0$ and $\mu_*(g) = \phi \ (\alev)$.
In \S\ref{sec:class:degenerate:green:function},
we will check that the birational class $\Br$
is a class of degeneration
(cf. Proposition~\ref{prop:Br:function:sp:pairing}).
Note that if $\dim X = 2$, then $\Br$ is nothing more than
$C^{\infty}(X(\CC), \RR)$.

The reason why we need to introduce the birational class
comes from the following observation.
Let $(E, h)$ be a Hermitian vector bundle of rank $2$ on $X$, and
$L$ a rank $1$ saturated subsheaf of $E$.
Then, $L_{\CC}$ is not necessarily line subbundle of $E_{\CC}$
if $\dim X_{\QQ} \geq 2$, so that
the metric $h_L$ of $L$ induced by $h$ is not 
necessarily $C^{\infty}$ over $X(\CC)$.
We can however see that
$\acherncl_1(L, h_L)$ gives rise to an element of
$\aBrChow^1(X)$.
Thus, we may consider
\[
\adeg \left( \acherncl_1(H, k)^{d-1} \cdot \acherncl_1(L, h_L) \cdot
(\acherncl_1(E, h) - \acherncl_1(L, h_L)) \right).
\]
We would like to compare the above with
$\adeg \left(  \acherncl_1(H, k)^{d-1} \cdot \acherncl_2(E, h) \right)$.
If we take the geometric case into account,
one can guess
\[
\adeg \left(  \acherncl_1(H, k)^{d-1} \cdot \acherncl_2(E, h) \right)
\geq \adeg \left( \acherncl_1(H, k)^{d-1} \cdot \acherncl_1(L, h_L) \cdot
(\acherncl_1(E, h) - \acherncl_1(L, h_L)) \right).
\]
Actually, we will prove the above inequality in 
\S\ref{sec:comp:intersection:numbers:via:birational:morphism}.
Using this together with the result of \cite{MoBGI},
we have the following Bogomolov's instability theorem for
rank $2$ vector bundles in the arithmetic case.

\begin{Theorem}[cf. Theorem~\ref{them:B:instability:rank:2}]
\label{them:B:instability:rank:2:intro}
If $\adeg\left( \acherncl_1(H, k)^{d-1} \cdot 
\left( 4 \acherncl_2(E, h) - \acherncl_1(E, h)^2 \right)
\right) < 0$, 
then there is a rank $1$ saturated subsheaf $L$ of $E$
such that $L^{\otimes 2} \otimes \det(E)^{-1}$ 
has positive degree on each connected component of $X(\CC)$
with respect to $H_{\CC}$, and that
\[
\adeg\left( \acherncl_1(H, k)^{d-1} \cdot 
\left( 2 \acherncl_1(L, h_L) - \acherncl_1(E, h) \right)^2
\right) > 0.
\]
\end{Theorem}

Finally, we would like to thank Prof. Bost who kindly
sent his paper \cite{Bo} after writing the first version of this note.
Inspired with his paper, we could simplify the description of this note.

\renewcommand{\theTheorem}{\arabic{section}.\arabic{Theorem}}
\renewcommand{\theClaim}{\arabic{section}.\arabic{Theorem}.\arabic{Claim}}
\renewcommand{\theequation}{\arabic{section}.\arabic{Theorem}.\arabic{Claim}}

\section{Class of degeneration of Green currents}
\label{sec:class:degenerate:green:function}
Let $X$ be a complex manifold.
We denote by $C^{\infty}(X)$ (resp. $C^{\infty}(X, \RR)$)
the set of all complex (resp. real) valued $C^{\infty}$ functions.
For a non-negative integer $k$ and a real number $r$ with 
$1 \leq r < \infty$,
we denote by $L^{r}_{k,loc}(X)$ (resp. $L^{r}_{k,loc}(X,\RR)$)
the set of all complex (resp. real) valued
functions on $X$ which locally have all weak partial derivatives up to order $k$ 
in $L^r$.
Let $\alpha$ be a form of type $(p,q)$ on $X$.
We say $\alpha$ is a locally $L^r_k$-form if
all coefficients of $\alpha$ in terms of
local coordinates belong to $L^r_{k,loc}$.
We denote the set of all (resp. real) locally 
$L^r_{k}$-forms of type $(p,q)$ on $X$
by $L^r_{k,loc}(\Omega_{X}^{p,q})$ 
(resp. $L^r_{k,loc}(\Omega_{X}^{p,q}, \RR)$).
Let us begin with the following lemma.

\begin{Lemma}
\label{lem:elem:prop:fun:sp:pairing}
\begin{enumerate}
\renewcommand{\labelenumi}{(\roman{enumi})}
\item
If $\phi \in L^r_{1, loc}(\Omega_X^{p,p})$ and
$\psi \in C^{\infty}(\Omega_X^{q,q})$, then
$\left[ \phi \wedge \partial\bar{\partial}(\psi) \right] + 
\left[ \partial(\phi) \wedge \bar{\partial}(\psi) \right]
\in \Image(\partial)$ and
$\left[ \phi \wedge \partial\bar{\partial}(\psi) \right] + 
\left[ \partial(\psi) \wedge \bar{\partial}(\phi) \right]
\in \Image(\bar{\partial})$ as currents.

\item
Let $r$ and $r'$ be real numbers with
$1 \leq r, r' < \infty$ and $1/r + 1/r' = 1$.
If $\phi \in L^r_{1,loc}(\Omega_X^{p,p})$ and $\psi \in L^{r'}_{1,loc}(\Omega_X^{q,q})$, 
then
$[\partial(\phi) \wedge \bar{\partial}(\psi)] =
[\partial(\psi) \wedge \bar{\partial}(\phi)]$
modulo $\Image(\partial) +  \Image(\bar{\partial})$.

\item
If $X$ is an $n$-dimensional connected compact complex manifold with
a fundamental form $\Phi$, then
\[
\sqrt{-1}\int_X 
\partial(\phi) \wedge \bar{\partial}(\phi) \wedge \Phi^{d-1} \geq 0
\]
for all $\phi \in L^2_{1}(X, \RR)$. Moreover, the equality holds if and only if
$\phi$ is a constant almost everywhere.
\end{enumerate}
\end{Lemma}

\Proof
(i) First of all, $\phi \wedge \bar{\partial}(\psi) \in L^r_{1, loc}(\Omega^{p+q,p+q+1}_X)$.
Thus, 
\[
\partial(\phi \wedge \bar{\partial}(\psi)) = \phi \wedge \partial\bar{\partial}(\psi) + 
\partial(\phi) \wedge \bar{\partial}(\psi).
\]
Hence, we get the first assertion. In the same way, 
\[
\bar{\partial}(\phi \wedge \partial(\psi)) = \phi \wedge \bar{\partial}\partial(\psi) + 
\bar{\partial}(\phi) \wedge \partial(\psi),
\]
which shows us the second assertion.

\medskip
(ii)
It is sufficient to see that
\[
\partial[\phi\wedge\bar{\partial}(\psi)] + \bar{\partial}[\phi\wedge\partial(\psi)]
= [\partial(\phi) \wedge \bar{\partial}(\psi)] -
[\partial(\psi) \wedge \bar{\partial}(\phi)]
\]
as currents.
This is a local question, so that we may assume that
$\phi$ and $\psi$ can be written in terms of a local coordinate and that all
coefficients of $\phi$ and $\psi$ belong to $L^r_1(X)$ and $L^{r'}_1(X)$
respectively.
Since $C^{\infty}(X) \cap L^{r'}_1(X)$ is dense in $L^{r'}_1(X)$,
there is a sequence $\{ \psi_n \}$ in $C^{\infty}(\Omega_X^{q,q})$
such that all coefficients of $\psi_n$ belong to $L^{r'}_1(X)$ and
they converge to the coefficients of $\psi$ in $L^{r'}_1(X)$.
Then,
\[
\begin{cases}
\lim_{n \to \infty}[\phi \partial(\psi_n)] = [\phi \partial(\psi)], \\
\lim_{n \to \infty}[\phi \bar{\partial}(\psi_n)] = 
[\phi \bar{\partial}(\psi)], \\
\lim_{n \to \infty}[\partial(\phi) \wedge \bar{\partial}(\psi_n)] =
[\partial(\phi) \wedge \bar{\partial}(\psi)], \\
\lim_{n \to \infty}[\partial(\psi_n) \wedge \bar{\partial}(\phi)] =
[\partial(\psi) \wedge \bar{\partial}(\phi)]
\end{cases}
\]
as currents. Here note that
if $T$ is a current, $\{ T_n \}$ is a sequence of currents, and
$\lim_{n \to \infty} T_n = T$ as currents, then
$\lim_{n \to \infty} \partial(T_n) = \partial(T)$ and
$\lim_{n \to \infty} \bar{\partial}(T_n) = \bar{\partial}(T)$
as currents. On the other hand, by virtue of the proof of (i),
\[
\partial(\phi\wedge\bar{\partial}(\psi_n)) + \bar{\partial}(\phi\wedge\partial(\psi_n))
= (\partial(\phi) \wedge \bar{\partial}(\psi_n)) -
(\partial(\psi_n) \wedge \bar{\partial}(\phi))
\]
for all $n$.
Thus, we get (ii).

\medskip
(iii)
Let $x$ be an arbitrary point of $X$, and $\theta_1, \ldots, \theta_n$
a local orthogonal frame of the holomorphic cotangent bundle
$\Omega^1_X$ around $x$ with respect to $\Phi$ such that
$\Phi = \sqrt{-1} \sum_i \theta_i \wedge \theta_i$.
If we set $\partial(\phi) = \sum_i a_i \theta_i$ around $x$, then
$\bar{\partial}(\phi) = \sum_i \bar{a}_i \bar{\theta}_i$.
Thus,
\[
\sqrt{-1} \partial(\phi) \wedge \bar{\partial}(\phi) \wedge \Phi^{n-1}
= \sum_{i=1}^n |a_i|^2 \Phi^n
\]
around $x$. This means that
\[
\sqrt{-1}\partial(\phi) \wedge \bar{\partial}(\phi)
\wedge \Phi^{n-1}
\]
is non-negative on $X$.
Therefore, we get the first assertion.

Next we assume the equality.
Then, by the proof of the inequality, we can see that
$\partial(\phi) = \bar{\partial}(\phi) = 0 \ (\alev)$, i.e.,
$d(\phi) = 0 \ (\alev)$.
Thus, $\phi$ is a constant almost everywhere.
\QED

An abelian group $B$ is called 
{\em a class of degeneration of Green currents for codimension $p$ cycles
in $L^r$} 
(or simply a class of degeneration)
if $C^{\infty}(\Omega_X^{p-1,p-1}, \RR) \subseteq B 
\subseteq L^{r}_{1,loc}(\Omega_X^{p-1,p-1}, \RR)$.
For example, $C^{\infty}(\Omega_X^{p-1,p-1}, \RR)$ and 
$L^r_{k, loc}(\Omega_X^{p-1,p-1}, \RR)$ ($k \geq 1$) are
classes of degeneration in $L^r$.

\bigskip
Let us consider a non-trivial example of 
class of degeneration of Green functions.
Let $\mu : Y \to X$ be a proper bimeromorphic morphism
of complex manifolds, $U$ the maximal open set of $X$ with
$\mu^{-1}(U) \overset{\sim}{\longrightarrow} U$, and
$\omega$ a form on $Y$. We define the form $\mu_*(\omega)$ on $X$
to be
\[
\mu_*(\omega)(x) =
\begin{cases}
\omega(\mu^{-1}(x)) & \text{if $x \in U$} \\
0                   & \text{if $x \not\in U$}.
\end{cases}
\]
Note that if $\omega$ is locally integrable, then
$\mu_*(\omega)$ is also locally integrable and
$\mu_*([\omega]) = [\mu_*(\omega)]$ as currents.
It is easy to see that $\mu_*(\omega_1 \wedge \omega_2) = 
\mu_*(\omega_1) \wedge \mu_*(\omega_2)$.

Let $D$ be a divisor on $X$.
A locally integrable function $g$ on $X$ is called
{\em a Green function for $D$} if
$g$ is $C^{\infty}$ over $X \setminus \Supp(D)$ and
$dd^c([g]) + \delta_D$ is represented by a $C^{\infty}$-form.
It is easy to see that if $g$ is a Green function for some
divisor, then for any points $x \in X$, there
are a meromorphic function $f$ around $x$ and
a $C^{\infty}$-function $\psi$ around $x$ with $g = \log |f| + \psi$.

Here we consider the following space $\Br(X)$.
A locally integrable function $\phi$ on $X$ belongs to $\Br(X)$
if and only if there are a proper bimeromorphic morphism
$\mu : Y \to X$ of complex manifolds, a divisor $D$ on $Y$, and
a Green function $g$ for $D$ such that
$\mu_*(D) = 0$ and $\phi = \mu_*(g)\ (\alev)$. 

\begin{Proposition}
\label{prop:Br:function:sp:pairing}
The space $\Br(X)$ is a class of degeneration of Green functions in $L^2$.
Moreover, the following properties are satisfied.
\begin{enumerate}
\renewcommand{\labelenumi}{(\arabic{enumi})}
\item
For all $\phi \in \Br(X)$, 
the differentials $\partial\bar{\partial}([\phi])$ in the sense of currents
are represented by locally integrable
forms. \rom{(}By abuse of notation,
representatives of $\partial\bar{\partial}([\phi])$ are denoted by 
$\partial\bar{\partial}(\phi)$.\rom{)}

\item
$\phi \wedge \partial\bar{\partial}(\psi)$ is
locally integrable forms for any $\phi, \psi \in \Br(X)$, and,
as currents,
$\left[ \phi \wedge \partial\bar{\partial}(\psi) \right] + 
\left[ \partial(\phi) \wedge \bar{\partial}(\psi) \right]
\in \Image(\partial)$.
\end{enumerate}
\end{Proposition}

\Proof
Obviously, $C^{\infty}(X, \RR) \subseteq \Br(X) \subseteq L^{2}_{loc}(X, \RR)$
because a $C^{\infty}$-function is a Green function for the zero divisor, 
and a Green function is locally square integrable.
First, let us check that $\Br(X)$ is an abelian group.
Choose arbitrary elements $\phi_1, \phi_2 \in \Br(X)$.
Then, we can easily find a proper bimeromorphic morphism
$\mu : Y \to X$ of complex manifolds, 
divisors $D_1, D_2$ on $Y$, and
Green functions $g_1$ for $D_1$ and $g_2$ for $D_2$ such that
$\mu_*(D_1) = \mu_*(D_2) = 0$,
$\phi_1 = \mu_*(g_1)\ (\alev)$, and
$\phi_2 = \mu_*(g_2)\ (\alev)$.
Thus,
$g_1 - g_2$
is a Green function for $D_1 - D_2$, and
$\mu_* (g_1 - g_2) = \phi_1 - \phi_2 \ (\alev)$. Hence,
$\phi_1 - \phi_2 \in \Br(X)$.
This shows us that $\Br(X)$ is an abelian group.

In order to check another properties, 
we need to prepare two lemmas.

\begin{Lemma}
\label{lem:eq:two:integrals:on:X}
Let $X$ be an $n$-dimensional complex manifold,
$D$ a divisor on $X$, and $g$ a Green function for $D$. Then,
\[
\int_X g d(\omega) = - \int_X d(g) \wedge \omega
\]
for all $\omega \in C^{\infty}_c(X, \Omega^{2n-1})$, i.e.,
$\omega$ is a compactly supported $(2n-1)$-form on $X$.
In other words, $d([g]) = [d(g)]$.
Note that $d(g)$ is a logarithmic form on $X$, so that
$d(g)$ is locally integrable.
\end{Lemma}

\Proof
Let $\mu : Y \to X$ be a proper bimeromorphic morphism such that
$\mu^{-1}(\Supp(D))$ is a normal crossing divisor.
Then,
\[
\int_X g d(\omega) =
\int_Y \mu^*(g) d(\mu^*(\omega))
\quad\text{and}\quad
\int_X d(g) \wedge \omega =
\int_Y d(\mu^*(g)) \wedge \mu^*(\omega).
\]
Thus, we may assume that $\Supp(D)$ is a normal crossing divisor.

Let $\{ U_{\alpha} \}_{\alpha \in A}$ be a locally finite open covering of $X$ 
such that each $U_{\alpha}$ is isomorphic to a bounded open set of $\CC^n$.
Let $\sum_{\alpha \in A} \phi_{\alpha} = 1$ be a partition of 
unity subordinate to
$\{ U_{\alpha} \}_{\alpha \in A}$. If
\[
\int_X g d(\phi_{\alpha} \omega) = 
- \int_X d(g) \wedge \phi_{\alpha} \omega
\]
for all $\alpha \in A$, then
{\allowdisplaybreaks
\begin{align*}
-\int_X d(g) \wedge \omega & 
= \sum_{\alpha \in A} -\int_X d(g) \wedge \phi_{\alpha} \omega \\
& = \sum_{\alpha \in A} \int_X g d(\phi_{\alpha} \omega) \\
& = \sum_{\alpha \in A} \int_X g \left( d(\phi_{\alpha}) \omega + \phi_{\alpha} d(\omega)
\right) \\
& =  \int_X g \left( 
d \left(  \sum_{\alpha \in A} \phi_{\alpha}\right) \omega +  
\sum_{\alpha \in A} \phi_{\alpha} d(\omega)
\right) \\
& = \int_X g d(\omega).
\end{align*}}
Thus, in order to complete our lemma,
it is sufficient to show the following sublemma.
\QED

\begin{Sublemma}
\label{sublem:eq:two:integrals:on:Cn}
Let $(z_1, \ldots, z_n)$ be a coordinate of $\CC^n$, and
$a_1, \ldots, a_n$ real numbers.
Then, for any $\omega \in C^{\infty}_c(\CC^n, \Omega^{2n-1})$,
\[
\int_{\CC^n} \log \left( |z_1|^{a_1} \cdots |z_n|^{a_n} \right) d(\omega) =
- \int_{\CC^n} d \left( \log \left( |z_1|^{a_1} \cdots |z_n|^{a_n} \right) \right) \wedge \omega.
\]
Note that
\[
d \left( \log \left( |z_1|^{a_1} \cdots |z_n|^{a_n} \right) \right) =
\sum_{i=1}^n \frac{a_i}{2} \left( \frac{dz_i}{z_i} + \frac{d \bar{z}_i}{\bar{z}_i} \right)
\]
and it is a $L^1$-form.
\end{Sublemma}

\Proof
Since $\log \left( |z_1|^{a_1} \cdots |z_n|^{a_n} \right) 
= a_1 \log|z_1| + \cdots + a_n \log|z_n|$,
it is sufficient to see that
\[
\int_{\CC^n} \log |z_1| d(\omega) =
- \int_{\CC^n} d \left( \log |z_1| \right) \wedge \omega.
\]
For $\epsilon > 0$, we set
\[
U_{\epsilon} = \{ (z_1, \ldots, z_n) \mid |z_1| \geq \epsilon \}\quad\text{and}\quad
D_{\epsilon} = \{ (z_1, \ldots, z_n) \mid |z_1| = \epsilon \}.
\]
Then, since $d(\log|z_1| \omega) = d(\log|z_1|) \omega + \log|z_1| d(\omega)$ over $U_{\epsilon}$,
by Stokes' formula,
\[
- \int_{D_{\epsilon}} \log|z_1| \omega = \int_{U_{\epsilon}} d(\log|z_1|) \wedge \omega +
\int_{U_{\epsilon}} \log|z_1| d(\omega).
\]
Moreover,
\[
\lim_{\epsilon \downarrow 0} \int_{U_{\epsilon}} d(\log|z_1|) \wedge \omega =
\int_{\CC^n} d(\log|z_1|) \wedge \omega
\quad\text{and}\quad
\lim_{\epsilon \downarrow 0} \int_{U_{\epsilon}} \log|z_1| d(\omega) =
\int_{\CC^n} \log|z_1| d(\omega).
\]
Thus, it is sufficient to show that
\[
\lim_{\epsilon \downarrow 0} \int_{D_{\epsilon}} \log|z_1| \omega = 0.
\]
Let us choose a sufficiently large number $M$ such that
$\operatorname{supp}(\omega) \subset \Delta_M^n$,
where $\Delta_M = \{ z \in \CC \mid |z| \leq M \}$.
Then, if we set $S_{\epsilon}^1 = \{ z \in \CC \mid |z| = \epsilon \}$, we have
\[
\int_{D_{\epsilon}} \log|z_1| \omega = \int_{S_{\epsilon}^1 \times \Delta_M^{n-1}} \log|z_1| \omega.
\]
Here we set
\begin{multline*}
\omega = \sum_{i=1}^n \left\{ \alpha_i 
(dx_1 \wedge dy_1) \wedge \cdots \wedge (\widehat{dx_i} \wedge dy_i) \wedge \cdots \wedge (dx_n \wedge dy_n) 
\right.\\
\left.
+ \beta_i
(dx_1 \wedge dy_1) \wedge \cdots \wedge (dx_i \wedge \widehat{dy_i}) \wedge \cdots \wedge (dx_n \wedge dy_n)
\right\},
\end{multline*}
where $z_i = x_i + \sqrt{-1}y_i$.
Then, since $dx_1 \wedge dy_1 = 0$ on $S^1_{\epsilon}$,
\begin{align*}
\int_{S_{\epsilon}^1 \times \Delta_M^{n-1}} \log|z_1| \omega & = 
\log(\epsilon) \int_{S_{\epsilon}^1 \times \Delta_M^{n-1}}
(\alpha_1 dy_1 + \beta_1 dx_1 ) \wedge (dx_2 \wedge dy_2) \wedge \cdots \wedge (dx_n \wedge dy_n) \\
& = \epsilon \log(\epsilon)
\int_{0}^{2\pi} \int_{\Delta_M^{n-1}} (\alpha_1^{\epsilon} \cos(\theta) - \beta_1^{\epsilon} \sin(\theta))
d\theta dx_2 dy_2 \cdots dx_n dy_n,
\end{align*}
where $\alpha_1^{\epsilon}(\theta, z_2, \cdots, z_n)  = \alpha_1(\epsilon e^{i \theta}, z_2, \cdots, z_n)$ and
$\beta_1^{\epsilon}(\theta, z_2, \cdots, z_n) = \beta_1(\epsilon e^{i \theta}, z_2, \cdots, z_n)$.
Thus, we get our lemma because 
\[
\lim_{\epsilon \downarrow 0} \epsilon \log(\epsilon) = 0
\quad\text{and}\quad
| \alpha_1^{\epsilon} \cos(\theta) - \beta_1^{\epsilon} \sin(\theta)|
\leq \Vert \alpha_1 \Vert_{\sup} + \Vert \beta_1 \Vert_{\sup}.
\]
\QED

Next, let us consider the following lemma.

\begin{Lemma}
\label{lem:smoothness:log:wedge:pullback}
Let $X$ be an $n$-dimensional complex manifold, and
\[
Y = X_N \overset{\pi_{N-1}}{\longrightarrow} X_{N-1} \overset{\pi_{N-2}}{\longrightarrow}
\cdots \overset{\pi_{0}}{\longrightarrow} X_0 = X
\]
a succession of blowing-ups along smooth and irreducible subvarieties of codimension
at least $2$, i.e., for each $0 \leq \alpha < N$, there is a smooth and irreducible
subvariety $C_{\alpha}$ on $X_{\alpha}$ such that
$\codim C_{\alpha} \geq 2$ and $\pi_{\alpha} : X_{\alpha + 1} \to X_{\alpha}$
is the blowing-up along $C_{\alpha}$.
Let $\Sigma$ be the exceptional set of 
$\pi = \pi_0 \cdots \pi_{N-1} : Y \to X$.
Let $D$ be a divisor on $Y$ with $\Supp(D) \subseteq \Sigma$, and
$g$ a Green function for $D$.
If $\Sigma$ is a normal crossing divisor, then
$d(g) \wedge \pi^*(\omega)$ is a $C^{\infty}$ form for any
$\omega \in C^{\infty}(X, \Omega^{n-1,n-1})$.
\end{Lemma}

\Proof
Let $y$ be an arbitrary point of $Y$, and
$(z_1, \ldots, z_n)$ a local coordinate of $Y$ at $y$
such that $z_1(y) = \cdots = z_n(y) = 0$ and
$\Sigma$ is given by $\{ z_1 \cdots z_a = 0 \}$ around $y$.
Then, $g$ can be written by a form 
\[
 g = e_1 \log |z_1|^2 + \cdots + e_a \log |z_a|^2 + 
(\text{$C^{\infty}$ function}).
\]
Then,
\[
d(g) = \sum_{i=1}^a e_i 
\left( \frac{dz_i}{z_i} + \frac{d\bar{z}_i}{\bar{z}_i} \right)
+ (\text{$C^{\infty}$ form}).
\]
Thus, it is sufficient to show that
\[
\frac{dz_i}{z_i} \wedge \pi^*(\omega)
\quad\text{and}\quad
\frac{d\bar{z}_i}{\bar{z}_i} \wedge \pi^*(\omega)
\]
are $C^{\infty}$ forms around $y$ for every $1 \leq i \leq a$.

We choose $0 \leq \alpha < N$ such that
$\{ z_i = 0 \}$ is an irreducible component of $\mu_{\alpha}^{-1}(C_{\alpha})$,
where $\mu_{\alpha} = \pi_{\alpha} \cdots \pi_{N-1} : Y \to X_{\alpha}$.
Moreover, we choose a local coordinate 
$(w_1, \ldots, w_n)$ of $X_{\alpha}$ at $\mu_{\alpha}(y)$
such that $w_1(\mu_{\alpha}(y)) = \cdots = w_n(\mu_{\alpha}(y)) = 0$ and
$C_{\alpha}$ is given by an equation $w_1 = \cdots = w_b = 0$.
Then, $b \geq 2$ because $\codim C_{\alpha} \geq 2$.
We set $\phi_i = \mu_{\alpha}^*(w_i)$ for $i = 1, \ldots, n$.
By our choice of $x_i$'s and $w_j$'s,
the ideal generated by $\phi_1, \ldots, \phi_b$ is contained in
the ideal generated by $z_i$.
Thus, there are holomorphic functions $f_1, \ldots, f_b$ around $y$
with $\phi_1 = z_i f_1, \ldots, \phi_b = z_i f_b$.
Here we set
\[
(\pi_{0} \cdots \pi_{\alpha-1})^*(\omega) =
\sum_{s,t} \omega_{st} \left( dw_1 \wedge \cdots \wedge \widehat{dw_s} \wedge \cdots \wedge dw_n \right)
\wedge \left( d\bar{w}_1 \wedge \cdots \wedge \widehat{d\bar{w}_t} \wedge \cdots \wedge d\bar{w}_n \right)
\]
around $\mu_{\alpha}(y)$.
Then,
\[
\pi^*(\omega) = 
\sum_{s,t} \mu_{\alpha}^*(\omega_{st}) 
\left( d\phi_1 \wedge \cdots \wedge \widehat{d\phi_s} \wedge \cdots \wedge d\phi_n \right)
\wedge \left( d\bar{\phi}_1 \wedge \cdots \wedge \widehat{d\bar{\phi}_t} \wedge \cdots \wedge d\bar{\phi}_n \right)
\]
Since $b \geq 2$, for each $s$,
there is $s'$ with $1 \leq s' \leq b$ and $s' \not= s$.
Then,
\[
\frac{dz_i}{z_i} \wedge d \phi_{s'} = 
\frac{dz_i}{z_i} \wedge ((dz_i) f_{s'} + z_i df_{s'}) =
dz_i \wedge df_{s'},
\]
which shows us that
\[
\frac{dz_i}{z_i} \wedge d\phi_1 \wedge \cdots \wedge \widehat{d\phi_s} \wedge \cdots \wedge d\phi_n
\]
is a holomorphic form for all $s$. In the same way, we can see that
\[
\frac{d\bar{z}_i}{\bar{z}_i} \wedge 
d\bar{\phi}_1 \wedge \cdots \wedge \widehat{d\bar{\phi}_t} \wedge \cdots \wedge d\bar{\phi}_n
\]
is an anti-holomorphic form for each $t$.
Thus, we get our lemma.
\QED

\renewcommand{\theClaim}{\arabic{section}.\arabic{Theorem}}

Let us go back to the proof of Proposition~\ref{prop:Br:function:sp:pairing}.
Let us pick up arbitrary $\phi, \psi \in \Br(X)$.
Choose a proper bimeromorphic morphism
$\mu : Y \to X$ of complex manifolds, divisors $D$ and $E$ on $Y$, and
Green functions $g$ for $D$ and $f$ for $E$ such that
$\mu_*(D) = \mu_*(E) = 0$, $\phi = \mu_*(g)\ (\alev)$, and
$\psi = \mu_*(f)\ (\alev)$.
Changing a model, if necessarily,
we may assume that $\mu : Y \to X$ can be obtained by
a succession of blowing-ups along smooth and irreducible subvarieties of codimension
at least $2$, and that the exceptional set is a divisor
with only normal crossings.
Let $\omega$ be a $C^{\infty}$-form on $Y$ with
\[
dd^c([f]) + \delta_E =
\frac{-1}{2 \pi \sqrt{-1}} \partial\bar{\partial}([f]) + \delta_E = [\omega].
\]

\addtocounter{Theorem}{1}
\begin{Claim}
\label{claim:prop:Br:function:sp:pairing:1}
$d([\phi]) = [\mu_*(d(g))]$ and 
$\partial\bar{\partial}([\psi]) = -2 \pi \sqrt{-1} [\mu_*(\omega)]$.
In particular, $d([\phi])$ and 
$\partial\bar{\partial}([\psi])$ 
are represented by locally integrable forms.
\end{Claim}

By virtue of Lemma~\ref{lem:eq:two:integrals:on:X},
$d([g]) = [d(g)]$. Thus, 
\[
d([\phi]) = d(\mu_*[g]) = \mu_* d([g]) = \mu_*[d(g)] = [\mu_*(d(g))].
\]
Further, 
$\mu_* \partial\bar{\partial}([f]) = -2 \pi \sqrt{-1} \mu_* [\omega]$
because $\mu_*(\delta_E) = 0$.
Thus, 
\[
\partial\bar{\partial}([\psi]) = \partial\bar{\partial}(\mu_* [f])  =
\mu_* \partial\bar{\partial}([f]) = -2 \pi \sqrt{-1} \mu_* [\omega]
= -2 \pi \sqrt{-1} [\mu_*(\omega)].
\]

\addtocounter{Theorem}{1}
\begin{Claim}
\label{claim:prop:Br:function:sp:pairing:2}
$\phi \partial\bar{\partial}(\psi)$,
$\partial(\phi) \wedge \bar{\partial}(\psi)$, and
$\phi \bar{\partial}(\psi)$ are locally integrable.
Note that the local integrability of
$\partial(\phi) \wedge \bar{\partial}(\phi)$ implies
that $\phi \in L^2_{1, loc}(X, \RR)$.
\end{Claim}

Using the equation 
$\partial\bar{\partial}(\psi) = -2 \pi \sqrt{-1} \mu_*(\omega)$,
we have
\[
\phi \partial\bar{\partial}(\psi) = -2 \pi \sqrt{-1} \mu_*(g) \mu_*(\omega) =
-2 \pi \sqrt{-1} \mu_*(g \omega).
\]
Here, since $g \omega$ is locally integrable,
so is $\mu_*(g \omega)$.
In order to see that $\partial(\phi) \wedge \bar{\partial}(\psi)$ is
locally integrable, it is sufficient to see that
$\partial(\phi) \wedge \bar{\partial}(\psi) \wedge \lambda$ is integrable
for all $\lambda \in C_c^{\infty}(\Omega_X^{n-1, n-1})$.
Since 
\[
\int_X \partial(\phi) \wedge \bar{\partial}(\psi) \wedge \lambda =
\int_Y \partial(g) \wedge \bar{\partial}(f) \wedge \mu^*(\lambda)
\]
and $\bar{\partial}(f) \wedge \mu^*(\lambda)$ is a $C^{\infty}$-form by
Lemma~\ref{lem:smoothness:log:wedge:pullback},
we can see that $\partial(\phi) \wedge \bar{\partial}(\psi) \wedge \lambda$ is integrable.
In the same way as above, using Lemma~\ref{lem:smoothness:log:wedge:pullback},
we can check that $\phi \bar{\partial}(\psi)$ are locally integrable.

\addtocounter{Theorem}{1}
\begin{Claim}
\label{claim:prop:Br:function:sp:pairing:3}
$\left[ \phi \partial\bar{\partial}(\psi) \right] + 
\left[ \partial(\phi) \wedge \bar{\partial}(\psi) \right] =
\partial[\phi \bar{\partial}(\psi)]$.
\end{Claim}

The above equation means that
\[
\int_X \phi \partial\bar{\partial}(\psi) \wedge \lambda + 
\int_X \partial(\phi) \wedge \bar{\partial}(\psi) \wedge \lambda =
\int_X \phi \bar{\partial}(\psi) \wedge \partial(\lambda)
\]
for all $\lambda \in C_c^{\infty}(\Omega_X^{n-1, n-1})$.
This is equivalent to say that
\[
-2\pi\sqrt{-1}\int_Y g \omega \wedge \mu^*(\lambda) + 
\int_X \partial(g) \wedge \bar{\partial}(f) \wedge \mu^*(\lambda) =
\int_X g \bar{\partial}(f) \wedge \partial(\mu^*(\lambda))
\]
for all $\lambda \in C_c^{\infty}(\Omega_X^{n-1, n-1})$.
We set $\eta = \bar{\partial}(f) \wedge \mu^*(\lambda)$.
Then, $\eta$ is $C^{\infty}$ by
Lemma~\ref{lem:smoothness:log:wedge:pullback} and
\[
d(\eta) = -2\pi\sqrt{-1}\omega \wedge \mu^*(\lambda) - 
\bar{\partial}(f) \wedge \mu^*(\partial(\lambda)).
\]
Thus, using Lemma~\ref{lem:eq:two:integrals:on:X},
\[
\int_Y \partial(g) \wedge \eta = \int_Y d(g) \wedge \eta =
- \int_Y g \wedge d(\eta) =
2\pi\sqrt{-1}\int_Y g \omega \wedge \mu^*(\lambda) + 
\int_Y g \bar{\partial}(f) \wedge \mu^*(\partial(\lambda)).
\]
Hence we get our claim.

\bigskip
Gathering Claim~\ref{claim:prop:Br:function:sp:pairing:1},
Claim~\ref{claim:prop:Br:function:sp:pairing:2} and 
Claim~\ref{claim:prop:Br:function:sp:pairing:3},
we can complete the proof of Proposition~\ref{prop:Br:function:sp:pairing}.
\QED

\renewcommand{\theTheorem}{\arabic{section}.\arabic{subsection}.\arabic{Theorem}}
\renewcommand{\theClaim}{\arabic{section}.\arabic{subsection}.\arabic{Theorem}.\arabic{Claim}}
\renewcommand{\theequation}{\arabic{section}.\arabic{subsection}.\arabic{Theorem}.\arabic{Claim}}

\section{Degenerate Green currents}

\subsection{$B$-Green currents and their star product}
\setcounter{Theorem}{0}
\label{subsec:green:fun:via:class}
Let $X$ be a complex manifold, and $B$ a 
class of degeneration of Green currents for codimension $p$ cycles.
Let $Z$ be a cycle of codimension $p$ on $X$.
A current (resp. locally integrable form) $g$ of type $(p-1,p-1)$ on $X$
is called a {\em $B$-Green current for $Z$} (resp. {\em $B$-Green form for $Z$}) 
if there are
a Green current (resp. Green form) $f$ for $Z$ and $\phi \in B$ with $g = f + [\phi]$
(resp. $g = f + \phi\ (\alev)$).
We denote $dd^c(g) + \delta_Z$ by $\omega(g)$.
For example, if $B = C^{\infty}(\Omega_X^{p-1,p-1}, \RR)$, then
a $B$-Green current is nothing more than an usual Green current.

We also fix a class $B'$ of degeneration of Green currents for codimension $q$ cycles.
We assume that
$B \subseteq L^r_{1,loc}(\Omega_X^{p-1,p-1}, \RR)$,
$B' \subseteq L^{r'}_{1,loc}(\Omega_X^{q-1,q-1}, \RR)$, and
$1/r + 1/r' = 1$.
Let $Z_1$ be a cycle of codimension $p$ on $X$, and
$Z_2$ a cycle of codimension $q$ on $X$.
Let $g_1$ be a $B$-Green current for $Z_1$,
and $g_2$ a $B'$-Green current for $Z_2$.
Let us choose a Green current $f_1$ for $Z_1$, 
a Green current $f_2$ for $Z_2$,
$\phi_1 \in B$, and $\phi_2 \in B'$ such that
$g_1 = f_1 + \phi_1$ and $g_2 = f_2 + \phi_2$.
We suppose that $Z_1$ and $Z_2$ intersect properly.
We would like to define the star product $g_1 * g_2$ of $g_1$ and $g_2$ to be
\[
g_1 * g_2 = f_1 * f_2 + [\omega(f_1) \wedge \phi_2] + [\phi_1 \wedge \omega(f_2)] -
\frac{\sqrt{-1}}{2\pi}[\partial(\phi_1) \wedge \bar{\partial}(\phi_2)]
\]
as an element of $D^{1,1}(X)/(\Image(\partial) + \Image(\bar{\partial}))$.
The following proposition says us that
the above product is well defined
and it is commutative modulo
$\Image(\partial) + \Image(\bar{\partial})$.

\begin{Proposition}
\label{prop:def:star:product:independence:commute}
\begin{enumerate}
\renewcommand{\labelenumi}{(\arabic{enumi})}
\item 
$g_1 * g_2$ modulo
$\Image(\partial) + \Image(\bar{\partial})$ is well defined, namely,
$f_1 * f_2 + [\omega(f_1) \wedge \phi_2] + [\phi_1 \wedge \omega(f_2)] -
\frac{\sqrt{-1}}{2\pi}[\partial(\phi_1) \wedge \bar{\partial}(\phi_2)]$ modulo
$\Image(\partial) + \Image(\bar{\partial})$
does not depend on the choices
of $f_1$, $f_2$, $\phi_1$ and $\phi_2$.

\item
$g_1 * g_2 = g_2 * g_1$ modulo
$\Image(\partial) + \Image(\bar{\partial})$.
\end{enumerate}
\end{Proposition}

\Proof
(1)
Let $g_1 = f'_1 + \phi'_1$ and $g_2 = f'_2 + \phi'_2$ be
another expressions of $g_1$ and $g_2$,
where $f'_1$ is a Green current for $Z_1$,
$f'_2$ is a Green current for $Z_2$, $\phi'_1 \in B$ and
$\phi'_2 \in B'$.
Then, there are smooth forms $\eta_1$ and $\eta_2$, and
currents $S_1$, $T_1$, $S_2$ and $T_2$ such that
\[
f'_1 = f_1 + \eta_1 + \partial(S_1) + \bar{\partial}(T_1)
\quad\text{and}\quad
f'_2 = f_2 + \eta_2 + \partial(S_2) + \bar{\partial}(T_2).
\]
Thus, we have
\[
\phi'_1 = \phi_1 - \eta_1 - \partial(S_1) - \bar{\partial}(T_1)
\quad\text{and}\quad
\phi'_2 = \phi_2 - \eta_2 - \partial(S_2) - \bar{\partial}(T_2).
\]
First of all, it is well known that
\addtocounter{Claim}{1}
\begin{equation}
\label{eqn:prop:def:star:product:independence:commute:1}
f'_1 * f'_2 = (f_1 + \eta_1) * (f_2 + \eta_2) \quad
\text{modulo $\Image(\partial) + \Image(\bar{\partial})$}.
\end{equation}
Moreover, since
$\omega(f'_1) = \omega(f_1 + \eta_1)$ and
$\omega(f'_1)$ is $\partial$ and $\bar{\partial}$-closed,
we can see that
\[
\omega(f'_1) \wedge \phi'_2 =
\omega(f_1 + \eta_1) \wedge (\phi_2 - \eta_2) -
\partial(\omega(f_1 + \eta_1) \wedge S_2) - 
\bar{\partial}(\omega(f_1 + \eta_1) \wedge T_2),
\]
which shows us that
\addtocounter{Claim}{1}
\begin{equation}
\label{eqn:prop:def:star:product:independence:commute:2}
\omega(f'_1) \wedge \phi'_2 =  \omega(f_1 + \eta_1) \wedge (\phi_2 - \eta_2)
\quad\text{modulo $\Image(\partial) + \Image(\bar{\partial})$}.
\end{equation}
In the same way,
\addtocounter{Claim}{1}
\begin{equation}
\label{eqn:prop:def:star:product:independence:commute:3}
\phi'_1 \wedge \omega(f'_2) =  (\phi_1 - \eta_1) \wedge \omega(f_2 + \eta_2)
\quad\text{modulo $\Image(\partial) + \Image(\bar{\partial})$}.
\end{equation}
Further, since
\[
\partial(\phi'_1) = \partial(\phi_1 - \eta_1) - \partial\bar{\partial}(T_1)
\quad\text{and}\quad
\bar{\partial}(\phi'_2) = \bar{\partial}(\phi_2 - \eta_2) - \bar{\partial}\partial(S_2),
\]
$\partial\bar{\partial}(T_1)$ (resp. $\bar{\partial}\partial(S_2)$) is
a $\partial$ and $\bar{\partial}$-closed locally $L^r_0$-form
(resp. $L^{r'}_0$-form).
Thus, we can see that
\[
\partial(\phi'_1) \wedge \bar{\partial}(\phi'_2) =
\partial(\phi_1 - \eta_1) \wedge \bar{\partial}(\phi_2 - \eta_2) +
\bar{\partial}(\partial\bar{\partial}(T_1) \wedge (\phi_2 - \eta_2)) +
\partial(\phi'_1 \wedge \partial\bar{\partial}(S_2)),
\]
which says us that
\addtocounter{Claim}{1}
\begin{equation}
\label{eqn:prop:def:star:product:independence:commute:4}
\partial(\phi'_1) \wedge \bar{\partial}(\phi'_2) =
\partial(\phi_1 - \eta_1) \wedge \bar{\partial}(\phi_2 - \eta_2)
\quad\text{modulo $\Image(\partial) + \Image(\bar{\partial})$}.
\end{equation}
Thus, gathering \eqref{eqn:prop:def:star:product:independence:commute:1},
\eqref{eqn:prop:def:star:product:independence:commute:2},
\eqref{eqn:prop:def:star:product:independence:commute:3}, and
\eqref{eqn:prop:def:star:product:independence:commute:4},
we obtain that
\[
f'_1 * f'_2 + [\omega(f'_1) \wedge \phi'_2] + [\phi'_1 \wedge \omega(f'_2)] -
\frac{\sqrt{-1}}{2\pi}[\partial(\phi'_1) \wedge \bar{\partial}(\phi'_2)]
\]
is equal to
\begin{multline*}
\Delta = (f_1 + \eta_1) * (f_2 + \eta_2) +  
[\omega(f_1 + \eta_1) \wedge (\phi_2 - \eta_2) ] + \\
[(\phi_1 - \eta_1) \wedge \omega(f_2 + \eta_2)] -
\frac{\sqrt{-1}}{2\pi}[\partial(\phi_1 - \eta_1) \wedge \bar{\partial}(\phi_2 - \eta_2)]
\end{multline*}
modulo $\Image(\partial) + \Image(\bar{\partial})$.
Moreover, by easy calculations, we can see that
\[
\Delta  - \left( f_1 * f_2 + [\omega(f_1) \wedge \phi_2] + [\phi_1 \wedge \omega(f_2)] -
\frac{\sqrt{-1}}{2\pi}[\partial(\phi_1) \wedge \bar{\partial}(\phi_2)] \right)
\]
is equal to
\begin{multline*}
\frac{\sqrt{-1}}{2\pi} \left([\phi_1 \wedge \partial\bar{\partial}(\eta_2)] 
+ [\partial(\phi_1) \wedge \bar{\partial}(\eta_2)]\right) +
(\eta_1 * f_2 - f_2 * \eta_1) + \\
\frac{\sqrt{-1}}{2\pi}\left([\partial\bar{\partial}(\eta_1) \wedge \phi_2] 
+ [\partial(\eta_1) \wedge \bar{\partial}(\phi_2)]\right) -
\frac{\sqrt{-1}}{2\pi}\left([\partial(\eta_1) \wedge \bar{\partial}(\eta_2)] 
+ [\eta_1 \wedge \partial\bar{\partial}\eta_2]\right),
\end{multline*}
which is elements of
$\Image(\partial) + \Image(\bar{\partial})$
by (i) of Lemma~\ref{lem:elem:prop:fun:sp:pairing} 
and \cite[Corollary~2.2.9]{GSArInt}.
Thus, we get (1).

\medskip
(2) It is well known that $f_1 * f_2 = f_2 * f_1$ modulo
$\Image(\partial) + \Image(\bar{\partial})$ 
(cf. \cite[Corollary~2.2.9]{GSArInt}).
Moreover, by (ii) of Lemma~\ref{lem:elem:prop:fun:sp:pairing}, 
$[\partial(\phi_1) \wedge  \bar{\partial}(\phi_2)]
= [\partial(\phi_2) \wedge \bar{\partial}(\phi_1)]$
modulo $\Image(\partial) + \Image(\bar{\partial})$.
Thus, we have (2).
\QED

\subsection{Birational Green function}
\setcounter{Theorem}{0}
Let $X$ be a complex manifold, and $D$ a divisor on $X$.
Let $\Br(X)$ be the class introduced in 
\S\ref{sec:class:degenerate:green:function}.
A $\Br(X)$-Green function for a divisor $D$ is specially  
called {\em a birational Green function for $D$}.
It is easy to see that $g$ is a birational Green function for $D$ if and
only if there are a proper bimeromorphic morphism
$\mu : Y \to X$ of complex manifolds, a divisor $D'$ on $Y$,
and a Green function $g'$ for $D'$ such that
$\mu_*(D') = D$ and $\mu_*(g') = g$.
A point at which $\omega(g)$ is not $C^{\infty}$ is called
{\em a singular point of $g$}.
In this sense, the pair $(D', g')$ is called 
{\em a resolution of singularities of $(D, g)$ by $\mu$}.

An idea of the birational Green functions arises from
the following observation.
Let $(E, h)$ a Hermitian vector bundle on $X$,
$L$ a rank-$1$ saturated subsheaf of $E$, and $s$ a rational
section of $L$.
If $\dim X \geq 2$, then $g = -\log h(s,s)$ is
not necessarily a Green function for $D = \zero(s)$
because $E/L$ is not locally free in general.
However, we can see that
$g$ is a birational Green function for $D$.
For, there are a proper bimeromorphic morphism
$\mu : Y \to X$ of complex manifolds, and
a rank-$1$ saturated subsheaf $L'$ of $\mu^*(E)$ such that
$\mu^*(E)/L'$ is locally free and
$\mu_*(L') = L$.
Since $L'$ is a sub-line bundle of $\mu^*(E)$,
$L'$ possesses a Hermitian metric $h'$ induced from
$\mu^*(h)$.
Let $s'$ be the rational section of $L'$ corresponding with $s$.
Then, by virtue of Poincar\'{e}-Lelong formula,
$g' = -\log h'(s', s')$ is a Green function for $D' = \zero(s')$.
Further, $D = \mu_*(D')$ and $g = \mu_*(g')$ by our construction.

By the above observation,
we have the following definition.
Let $X$ be a complex manifold,
and $E$ a torsion free sheaf on $X$.
We say $h$ is {\em a birationally Hermitian metric of $E$} if
there is a proper bimeromorphic morphism $\mu : Y \to X$ of
complex manifolds, 
a Hermitian vector bundle $(E', h')$ on
$Y$, and a Zariski open set $U$ of $X$ such that
$\codim(X \setminus U) \geq 2$, $\mu^{-1}(U) \simeq U$, 
$(E, h)$ is a Hermitian vector bundle on $U$, and
$\rest{(E', h')}{\mu^{-1}(U)} \simeq \rest{(E, h)}{U}$.
We say $(E', h')$ is a resolution of singularities 
of $(E,h)$ by $\mu$.
Then, we have the following proposition.

\begin{Proposition}
\label{prop:criterion:birat:herm}
Let $L$ be a line bundle on $X$ and
$h$ a generalized metric of $L$ over $X$.
Let $s$ be a non-zero rational section of $L$.
Then, $h$ is a birationally Hermitian metric if and only if
$-\log h(s,s)$ is a birational Green function for
$\zero(s)$.
\end{Proposition}

\Proof
First, we assume that $h$ is a birationally Hermitian metric.
Let $(L', h')$ be a resolution of singularities of
$(L, h)$ by $\mu : Y \to X$. Let $U$ be a Zariski open set
of $X$ as the above.
Here, we can find a rational section $s'$ of $L'$ 
corresponding with $s$ via isomorphism
$\rest{(L', h')}{\mu^{-1}(U)} \simeq \rest{(L, h)}{U}$.
Then,
$\mu_*(\zero(s')) = \zero(s)$, and
$-\log h'(s', s')$ is a Green function for $\zero(s')$.
Thus, $-\log h(s,s)$ is a birational Green function for
$\zero(s)$ because $-\log h(s,s) = \mu_*(-\log h'(s', s'))$.

Next, we assume that
$g = -\log h(s,s)$ is a birational Green function for
$D = \zero(s)$. Let $(D', g')$ be a resolution of
singularities of $(D, g)$ by $\mu : Y \to X$.
Let $s'$ be a non-zero rational section of $\OO_Y(D')$ with
$D' = \zero(s')$.
Let $U$ be a non-empty Zariski open set of $X$
with $\mu^{-1}(U) \simeq U$.
Since $\rest{D'}{\mu^{-1}(U)} = \rest{D}{U}$,
there is $u \in H^0(U, \OO_U)^{\times}$ with $s' = us$.
Thus, there is an isomorphism $\iota : \rest{\OO_Y(D')}{\mu^{-1}(U)} 
\overset{\sim}{\longrightarrow} \rest{\OO_X(D)}{U}$ with $\iota(s') = s$.
Since $g'$ is a Green function for $D'$, there is
a $C^{\infty}$ Hermitian metric $h'$ of $\OO_Y(D')$ with
$g' = -\log h'(s', s')$.
Then, we have $h'(s', s') = h(s, s)$ over
$U = \mu^{-1}(U)$ because $g = \mu_*(g')$.
Thus, $\iota$ gives rise to an isometry
$\rest{(\OO_Y(D'), h')}{\mu^{-1}(U)} \overset{\sim}{\longrightarrow}
\rest{(\OO_X(D), h)}{U}$.
\QED

\begin{Corollary}
\label{cor:birat:green:produce:birat:herm}
Let $g$ be a birational Green function for a divisor $D$.
Then, there is a birationally Hermitian metric $h$ of $\OO_X(D)$
such that $g = -\log h(s, s)$, where $s$ is a non-zero rational section
of $\OO_X(D)$ with $\zero(s) = D$.
\end{Corollary}

\Proof
Since $g$ is a locally integrable function, there is
a generalized metric $h$ of $\OO_X(D)$ with
$g = -\log h(s, s)$. Hence, by Proposition~\ref{prop:criterion:birat:herm},
$h$ must be a birationally Hermitian metric.
\QED

\section{Arithmetic $B$-cycles and their pairing}
\renewcommand{\theequation}{\arabic{section}.\arabic{subsection}.\arabic{Theorem}}

\subsection{Arithmetic $B$-cycles}
Let $X$ be an arithmetic variety, i.e.,
a quasi-projective integral scheme over $\ZZ$ with
the smooth generic fiber over $\ZZ$.
For each $p$,
let $\widehat{R}^p(X)$ be the subgroup of $\aCycle^p(X)$ generated
by the following elements:
\begin{enumerate}
\renewcommand{\labelenumi}{(\alph{enumi})}
\item 
$((f), - [\log |f|^2])$, 
where $f$ is a rational function on some
subvariety $Y$ of codimension $p-1$ and $[\log |f|^2]$ 
is the current defined by
\[
[\log |f|^2](\gamma) = 
        \int_{Y(\CC)} (\log |f|^2)\gamma.
\]

\item
$(0, \partial(\alpha) + \bar{\partial}(\beta))$,
where $\alpha \in D^{p-2, p-1}(X(\CC))$,
$\beta \in D^{p-1, p-2}(X(\CC))$.
\end{enumerate}

A pair $(Z, g)$ is
called {\em an arithmetic $D$-cycle on $X$ of codimension $p$}
if $Z$ is a cycle of codimension $p$ on $X$ and
$g \in D^{p-1,p-1}(X(\CC))$.
We denote the set of all arithmetic $D$-cycles of codimension $p$ on $X$
by $\aDCycle^p(X)$. Moreover, $\aDChow^p(X)$ is defined by
$\aDCycle^{p}(X)/\widehat{R}^p(X)$.

Let $B$ be a class of degeneration of Green currents for codimension $p$
cycles on $X(\CC)$.
A pair $(Z, g)$ is called {\em an arithmetic $B$-cycle of codimension $p$ on $X$}
if $Z \in \Cycle^p(X)$ and $g$ is a $B$-Green current
for $Z(\CC)$ on $X(\CC)$.
We denote the set of all arithmetic $B$-cycles of codimension $p$ on $X$
by $\aBCycle^p(X)$, and define
\[
\aBChow^p(X) = \aBCycle^{p}(X)/\widehat{R}^p(X).
\]
If $B = L^r_{k,loc}(\Omega^{p-1,p-1}_{X(\CC)}, \RR)$, then
$\aCycle^p_B(X)$ and $\aChow^p_B(X)$ are denoted by
$\aCycle^p_{L^r_k}(X)$ and $\aChow^p_{L^r_k}(X)$ respectively.
Note that for $x \in \aChow^p_B(X)$, there are $x_0 \in \aChow^p(X)$ and
$\phi \in B$ with $x = x_0 + a(\phi)$.

\bigskip
From now, we assume that $X$ is regular.
We would like to construct a pairing
\addtocounter{Theorem}{1}
\begin{equation}
\label{eqn:pair:cycle:dg}
\aChow_{L^r_1}^p(X) \otimes \aChow_{L^{r'}_1}^q(X) \to \aDChow^{p+q}(X)_{\QQ},
\end{equation}
where $1/r + 1/r' = 1$.
Roughly speaking, for $(Z_1, g_1) \in \aCycle^p_{L^r_1}(X)$ and
$(Z_2, g_2) \in \aCycle_{L^{r'}_1}^q(X)$,
$(Z_1, g_1) \cdot (Z_2, g_2)$ is defined by $(Z_1 \cdot Z_2, g_1 * g_2)$ if
$Z_1$ and $Z_2$ intersect properly.

To define it exactly, 
for $x \in \aChow^p_{L^r_1}(X)$ and $y \in \aChow^q_{L^{r'}_1}(X)$,
we choose
$x_0 \in \aChow^p(X)$, $y_0 \in \aChow^q(Y)$,
$\phi \in L^r_{k,loc}(\Omega^{p-1,p-1}_{X(\CC)}, \RR)$, and
$\psi \in L^{r'}_{k,loc}(\Omega^{q-1,q-1}_{X(\CC)}, \RR)$
with $x = x_0 + a(\phi)$ and $y = y_0 + a(\psi)$.
Then, we define $x \cdot y$ to be
\[
x \cdot y = x_0 \cdot y_0 +
a \left(
[\omega(x_0) \wedge \psi] + [\phi \wedge \omega(y_0)] -
\frac{\sqrt{-1}}{2\pi}[\partial(\phi) \wedge \bar{\partial}(\psi)] 
\right),
\]
where $x_0 \cdot y_0$ is the usual arithmetic intersection.
In the same way as in the proof of 
Proposition~\ref{prop:def:star:product:independence:commute},
we can see that the above definition 
does not depend on the choices
of $x_0$, $y_0$, $\phi$ and $\psi$.
Moreover, the pairing is commutative by virtue of
(ii) of Lemma~\ref{lem:elem:prop:fun:sp:pairing}.

\bigskip
Here we recall the scalar product
\addtocounter{Theorem}{1}
\begin{equation}
\label{eqn:scalar:product:D:cycle}
\aChow^p(X) \otimes \aDChow^q(X) \to \aDChow^{p+q}(X)_{\QQ}
\end{equation}
introduced in \cite[\S2.3]{KMRB}.
Let $x \in \aChow^p(X)$ and $y \in \aDChow^q(X)$.
We choose $y_0 \in \aChow^q(X)$ and $T \in D^{p-1,p-1}(X(\CC))$
with $y = y_0 + a(T)$.
Then, $x \cdot y$ is defined by
\[
x \cdot y = x \cdot y_0 + a(\omega(x) \wedge T).
\]
If $y \in \aChow^q_{L^r_k}(X)$, then $T$ can be represented by
$\phi \in L^r_{k,loc}(\Omega^{q-1,q-1}_{X(\CC)}, \RR)$.
Thus, $\omega(x) \wedge T \in L^r_{k,loc}(\Omega_{X(\CC)}^{p+q-1,p+q-1}, \RR)$.
This observation shows us that \eqref{eqn:scalar:product:D:cycle} induces
\addtocounter{Theorem}{1}
\begin{equation}
\label{eqn:scalar:product:Lrk:cycle}
\aChow^p(X) \otimes \aChow^q_{L^r_k}(X) \to \aChow^{p+q}_{L^r_k}(X)_{\QQ}
\end{equation}
Then, we have the following proposition.

\begin{Proposition}
\label{prop:assoc:pairing}
Let $r$ and $r'$ be real numbers with $1 \leq r, r' < \infty$ and
$1/r + 1/r' = 1$.
For $x \in \aChow^p(X)$, $y \in \aChow^q_{L^r_1}(X)$, and
$z \in \aChow^s_{L^{r'}_1}(X)$,
$x \cdot (y \cdot z) = (x \cdot y) \cdot z$.
\end{Proposition}

\Proof
We set $y = y_0 + a(\phi)$ and $z = z_0 + a(\psi)$, where
$y_0 \in \aChow^q(X)$, $z_0 \in \aChow^s(Y)$,
$\phi \in L^r_{k,loc}(\Omega^{q-1,q-1}_{X(\CC)}, \RR)$, and
$\psi \in L^{r'}_{k,loc}(\Omega^{s-1,s-1}_{X(\CC)}, \RR)$.
Then,
\begin{align*}
x \cdot (y \cdot z) & = x \cdot 
\left(y_0 \cdot z_0 + a \left(
[\omega(y_0) \wedge \psi] + [\phi \wedge \omega(z_0)] -
\frac{\sqrt{-1}}{2\pi}[\partial(\phi) \wedge \bar{\partial}(\psi)] 
\right) \right) \\
& = x \cdot (y_0 \cdot z_0) + a \left( \omega(x) \wedge \left(
[\omega(y_0) \wedge \psi] + [\phi \wedge \omega(z_0)] -
\frac{\sqrt{-1}}{2\pi}[\partial(\phi) \wedge \bar{\partial}(\psi)] \right) 
\right).
\end{align*}
On the other hand,
\begin{align*}
(x \cdot y) \cdot z & = (x \cdot y_0 + a(\omega(x) \wedge \phi)) \cdot z \\
& = (x \cdot y_0) \cdot z_0 + \\
& \qquad\qquad
a \left( [\omega(x \cdot y_0) \wedge \psi] + 
[\omega(x) \wedge \phi \wedge \omega(z_0)]  -
\frac{\sqrt{-1}}{2\pi}
[\partial(\omega(x) \wedge \phi) \wedge \bar{\partial}(\psi)] \right).
\end{align*}
Here $\omega(x \cdot y_0) = \omega(x) \wedge \omega(y_0)$ and
$\partial(\omega(x) \wedge \phi) = \omega(x) \wedge \partial(\phi)$.
Thus, we have our proposition.
\QED

\subsection{Intersection on singular varieties}
\setcounter{Theorem}{0}
Let $X$ be an arithmetic variety.
Let $\overline{E} = (E, h)$ be a Hermitian vector bundle on $X$.
In \cite[Theorem~4]{GSRR}, the operator
\addtocounter{Theorem}{1}
\begin{equation}
\label{eqn:operator:A:cycle}
\achernch(\overline{E}) : \aChow^*(X) \to \aChow^*(X)_{\QQ}
\end{equation}
is defined. Here we would like to extend the above operator in two ways.

The first one is
\addtocounter{Theorem}{1}
\begin{equation}
\label{eqn:operator:D:cycle}
\achernch(\overline{E}) : \aDChow^*(X) \to \aDChow^*(X)_{\QQ}.
\end{equation}
This is defined by
\[
\achernch(\overline{E}) \cdot (x_0 + a(T)) = \achernch(\overline{E}) \cdot x_0 +
a(\chernch(\overline{E}) \wedge T),
\]
where $x_0 \in \aChow^*(X)$ and $T \in \bigoplus_{p \geq 1} D^{p-1,p-1}(X(\CC))$.

To give the second one, 
let us fix 
an element $\phi$ of $\bigoplus_{p \geq 1} L^r_{1,loc}(\Omega_{X(\CC)}^{p-1,p-1}, \RR)$.
We define
\addtocounter{Theorem}{1}
\begin{equation}
\label{eqn:operator:B:cycle}
\achernch(\overline{E}) + a(\phi): \aChow^*_{L^{r'}_1}(X) \to \aDChow^*(X)_{\QQ}
\end{equation}
as follows, where $1/r + 1/r' = 1$.
Let $x \in \aChow^*_{L^{r'}_1}(X)$. We choose $x_0 \in \aChow^*(X)$ and
$\psi \in \bigoplus_{p \geq 1} L^{r'}_{1,loc}(\Omega_{X(\CC)}^{p-1,p-1}, \RR)$
with $x = x_0 + a(\psi)$.
Then,
\[
(\achernch(\overline{E}) + a(\phi)) \cdot x =
\achernch(\overline{E}) \cdot x_0 + a\left(
[\chernch(\overline{E}) \wedge \psi] + [\phi \wedge \omega(x_0)] -
\frac{\sqrt{-1}}{2\pi}[\partial(\phi) \wedge \bar{\partial}(\psi)] 
\right).
\]

In the same way as in Proposition~\ref{prop:assoc:pairing}, 
using \cite[4 of Theorem~4]{GSRR}, we have the following.

\begin{Proposition}
\label{prop:assoc:operator}
Let $\overline{E}$ and $\overline{F}$ be Hermitian vector bundles on $X$.
Then, we have the following.
\begin{enumerate}
\renewcommand{\labelenumi}{(\arabic{enumi})}
\item
$\achernch(\overline{E}) \cdot (\achernch(\overline{F}) \cdot x) =
\achernch(\overline{E} \otimes \overline{F}) \cdot x$
for $x \in \aDChow^*(X)$.

\item
Let $r$ and $r'$ be real numbers with $1 \leq r, r' < \infty$ and
$1/r + 1/r' = 1$. Then, for any 
$\phi \in \bigoplus_{p \geq 1} L^r_{1,loc}(\Omega_{X(\CC)}^{p-1,p-1}, \RR)$ and
$x \in \aChow^*_{L^{r'}_1}(X)$,
\[
\achernch(\overline{E}) \cdot \left( (\achernch(\overline{F}) + a(\phi)) \cdot x \right)
= \left( \achernch(\overline{E} \otimes \overline{F}) + 
a(\chernch(\overline{E}) \wedge \phi) \right) \cdot x.
\]
\end{enumerate}
\end{Proposition}

\subsection{Arithmetic $B$-Cartier divisor}
\setcounter{Theorem}{0}
Let $X$ be an arithmetic variety.
Let $\Rat_X$ be the sheaf of rational functions on $X$.
We denote $H^0(X, \Rat_X^{\times}/\OO_X^{\times})$ by $\Div(X)$.
An element of $\Div(X)$ is called {\em a Cartier divisor on $X$}.
For a Cartier divisor $D$ on $X$, we can assign the divisor $[D] \in \Cycle^1(X)$
in the natural way. This gives rise to a homomorphism 
\[
c_X : \Div(X) \to \Cycle^1(X).
\]
Note that $c_X$ is neither injective nor surjective
in general. The exact sequence
\[
1 \to \OO_X^{\times} \to \Rat_X^{\times} \to \Rat_X^{\times}/\OO_X^{\times} \to 1
\]
induces to a homomorphism $\Div(X) \to H^1(X, \OO_X^{\times})$.
For a Cartier divisor $D$ on $X$, the image of $D$ by
the above homomorphism induces a line bundle on $X$.
We denote this line bundle by $\OO_X(D)$.
{\em An arithmetic Cartier divisor on $X$} is a pair $(D, g)$
such that $D \in \Div(X)$ and $g$ is a Green function
for $D(\CC)$ on $X(\CC)$. The set of all arithmetic Cartier divisors on $X$
is denoted by $\aDiv(X)$, and
$\aPic(X)$ is defined by $\aDiv(X)/\hat{c}_X^{-1}(\widehat{R}^1(X))$,
where $\hat{c}_X$ is a natural homomorphism
$\aDiv(X) \to \aCycle^1(X)$.

Let us fix a class $B$ of degeneration of Green functions
on $X(\CC)$, namely, an abelian group with
$C^{\infty}(X, \RR) \subseteq B \subseteq L^r_{1,loc}(X, \RR)$.
A pair $(D, g)$ is called {\em an arithmetic $B$-Cartier divisor on $X$}
if $D \in \Div(X)$ and $g$ is a $B$-Green function
for $D(\CC)$ on $X(\CC)$.
We denote the set of all arithmetic $B$-Cartier divisors on $X$
by $\aBDiv(X)$, and define
\[
\aBPic(X) = \aBDiv(X)/\hat{c}_X^{-1}(\widehat{R}^1(X)).
\]
Note that if $X$ is regular, then
$\aBDiv(X) = \aBCycle^1(X)$ and $\aBPic(X) = \aBChow^1(X)$.

If $B = L^r_{k,loc}(\Omega_{X(\CC)}^{p-1,p-1}, \RR)$,
then $\aBDiv(X)$ and $\aBPic(X)$
are denoted by 
$\aDiv_{L^r_k}(X)$ and $\aPic_{L^r_k}(X)$.
Moreover,
If $B = \Br(X(\CC))$, then
$\aBCycle^1(X)$, $\aBDiv(X)$, $\aBChow^1(X)$, and $\aBPic(X)$
are denoted by 
$\aBrCycle^1(X)$, $\aBrDiv(X)$, $\aBrChow^1(X)$, and $\aBrPic(X)$.
An element of $\aBrCycle^1(X)$ (resp. $\aBrDiv(X)$)
is called a {\em birational arithmetic divisor}
(resp. a {\em birational arithmetic Cartier divisor}).

We can easily to see that \eqref{eqn:operator:D:cycle} and
\eqref{eqn:operator:B:cycle} induce
\addtocounter{Theorem}{1}
\begin{equation}
\label{eqn:operator:D:Pic}
\aPic(X) \otimes \aDChow^p(X) \to \aDChow^{p+1}(X).
\end{equation}
and
\addtocounter{Theorem}{1}
\begin{equation}
\label{eqn:pair:C:div:cycle}
\aPic_{L^r_1}(X) \otimes \aChow^p_{L^{r'}_1}(X) \to \aDChow^{p+1}(X)
\end{equation}
respectively, where $1/r + 1/r' = 1$.
Note that if $D$ and $Z$ intersect properly, then
\eqref{eqn:pair:C:div:cycle} is given by
$(D, g_D) \cdot (Z, g_Z) = (D \cdot Z, g_D * g_Z)$,
where $D \cdot Z$ is defined as follows.
Let $s$ be a rational section of $\OO_X(D)$ with $\zero(s) = D$, and
$Z = a_1 Z_1 + \cdots + a_n Z_n$ the irreducible decomposition
as cycles. Then, $\rest{s}{Z_i}$ gives rise to a rational section of
$\OO_{Z_i}(D)$, so that
we define $D \cdot Z$ to be
\[
a_1 \zero( \rest{s}{Z_1} ) + \cdots + a_n \zero( \rest{s}{Z_n} ).
\]
In particular, if $B \subseteq L^2_{1,loc}(X, \RR)$, then we have
a commutative pairing:
\addtocounter{Theorem}{1}
\begin{equation}
\label{eqn:pair:birat:div}
\aBPic(X) \otimes \aBPic(X) \to \aDChow^2(X).
\end{equation}

\begin{Remark}
\label{rem:intersection:Br}
Let $x_0 \in \aPic(X)$, $y_0 \in \aPic(X)$,
$\phi \in \Br(X(\CC))$, and $\psi \in \Br(X(\CC))$.
Then,
\[
(x_0 + a(\phi)) \cdot (y_0 + a(\psi)) = x_0 \cdot y_0 +
a(\omega(x_0)\psi + \phi \omega(y_0) + dd^c(\phi)\psi).
\]
Indeed, by Proposition~\ref{prop:Br:function:sp:pairing} and
(ii) of Lemma~\ref{lem:elem:prop:fun:sp:pairing}, we have
\[
[\partial\bar{\partial}(\phi) \psi] =
-[ \partial(\phi) \wedge \bar{\partial}(\psi)]
= - [ \partial(\phi) \wedge \bar{\partial}(\psi)]
\]
modulo $\Image(\partial) + \Image(\bar{\partial})$.
\end{Remark}

\renewcommand{\theTheorem}{\arabic{section}.\arabic{Theorem}}
\renewcommand{\theClaim}{\arabic{section}.\arabic{Theorem}.\arabic{Claim}}
\renewcommand{\theequation}{\arabic{section}.\arabic{Theorem}.\arabic{Claim}}
\section{Hodge index theorem}
\label{sec:hodge:index:thm}

Let $K$ be a number field, and $O_K$ the ring of integers of $K$.
Let $f : X \to \Spec(O_K)$ be a projective arithmetic variety
with the geometrically irreducible generic fiber.
Let $\overline{H} = (H, k)$ be an arithmetically ample
Hermitian line bundle on $X$, i.e.,
(1) $H$ is $f$-ample,
(2) the Chern form $c_1(H, k)$ is 
positive definite on the infinite fiber $X(\CC)$, and
(3) there is a positive integer $m_0$ such that,
for any integer $m \geq m_0$, $H^0(X, H^m)$
is generated by the set
$\left\{ s \in H^0(X, H^m) \mid 
\Vert s \Vert_{\sup} < 1 \right\}$.
Let us consider the pairing
\[
(\ \cdot \ )_{\overline{H}} : \aPic_{L^2_1}(X)_{\QQ} \times \aPic_{L^2_1}(X)_{\QQ}
\to \RR
\]
given by $(x \cdot y)_{\overline{H}} =
\adeg \left( \acherncl_1(H, k)^{d-1} \cdot (x \cdot y) \right)$,
where $d = \dim X_K$.
Further, we have the homomorphism
\[
\deg_K : \aPic_{L^2_1}(X)_{\QQ} \to \QQ
\]
given by $\deg_K((D, g)) = ( D_K \cdot c_1(H_K)^{d-1})$.

\begin{Theorem}
\label{thm:hodge:index:theorem}
If $x \in \aPic_{L^2_1}(X)_{\QQ}$ and $\deg_K(x) = 0$, then
$(x \cdot x)_{\overline{H}} \leq 0$. Moreover,
the equality holds if and only if there is 
$y \in \aChow^1(\Spec(O_K))_{\QQ}$ with $x = f^*(y)$
in $\aChow^1(X)_{\QQ}$.
\end{Theorem}

\Proof
Clearly, we may assume that $x \in \aPic_{L^2_1}(X)$.
There is $x_0 \in \aPic(X)$ such that
$z(x) = z(x_0)$ and $\omega(x_0)$ is harmonic with
respect to $c_1(H, k)$.
Consequently, $\omega(x_0) \wedge c_1(H, k)^{d-1} = 0$
because $\deg_K(x) = 0$.
Then, we can find $\phi \in L^2_{1,loc}(X, \RR)$ with
$x = x_0 + a(\phi)$.
Thus, since $x_0 \cdot a(\phi) = a(\phi \omega(x_0) )$ and
$a(\phi) \cdot a(\phi) = 
a\left(\frac{-\sqrt{-1}}{2\pi}\partial(\phi) \wedge 
\bar{\partial}(\phi)\right)$, we have
\begin{align*}
(x \cdot x)_{\overline{H}} & = (x_0 \cdot x_0)_{\overline{H}} 
+ 2(x_0 \cdot a(\phi))_{\overline{H}}
+ (a(\phi) \cdot a(\phi))_{\overline{H}} \\
& = (x_0 \cdot x_0)_{\overline{H}} 
- \frac{\sqrt{-1}}{4\pi} \int_{X(\CC)} \partial(\phi) \wedge \bar{\partial}(\phi) \wedge c_1(H, k)^{d-1}.
\end{align*}
On the other hand,
by (ii) of Lemma~\ref{lem:elem:prop:fun:sp:pairing},
\[
\sqrt{-1} \int_{X(\CC)} 
\partial(\phi) \wedge \bar{\partial}(\phi) \wedge c_1(H, k)^{d-1}
\geq 0.
\]
Thus, $(x \cdot x)_{\overline{H}} \leq (x_0 \cdot x_0)_{\overline{H}}$.
Moreover, by \cite[Theorem~1.1]{MoHI}, $(x_0 \cdot x_0)_{\overline{H}} \leq 0$.
Therefore, we get the first assertion.

Next we assume that the equality holds. Then,
\[
\sqrt{-1} \int_{X(\CC)} 
\partial(\phi) \wedge \bar{\partial}(\phi) \wedge c_1(H, k)^{d-1}
= 0.
\]
Thus, by the equality condition of
(ii) of Lemma~\ref{lem:elem:prop:fun:sp:pairing},
there is a function $c$ on $X(\CC)$ such that
$\phi = c \ (\alev)$ and
$c$ is constant on each connected component of $X(\CC)$.
Therefore,
\[
x = x_0 + a(c) \in \aPic(X).
\]
Hence, by \cite[Theorem~1.1]{MoHI},
then there are $(D, g) \in \aDiv(X)$ and
a positive integer $n$ such that $nx \sim (D, g)$,
$D$ is vertical with respect to $X \to \Spec(O_K)$ and
$g$ is constant on each connected component of $X(\CC)$.
Then,
\[
\adeg\left( \acherncl_1(H, k)^{d-1} \cdot (D, g)^2 \right) =
\sum_{P \in \Spec(O_K) \setminus \{ 0 \}} \deg(H_P^{d-1} \cdot D_P^2).
\]
Thus, by Zariski's lemma for integral scheme
(cf. Lemma~\ref{lem:int:on:fibers}),
there is a $\QQ$-divisor $T$ on $\Spec(O_K)$ with
$f^*(T) = D$ in $\Cycle^1(X)_{\QQ}$. 
Let $K(\CC)$ be the set of all embeddings of $K$ into $\CC$, and,
for each $\sigma \in K(\CC)$,
let $X_{\sigma} = X \otimes_K^{\sigma} \CC$ be the base extension in terms of
$\sigma$. Then,
$X(\CC) = \coprod_{\sigma} X_{\sigma}$ is nothing more than
the decomposition into connected component.
Let $g_{\sigma}$ be the value of $g$ on $X_{\sigma}$.
Then, $(D, g) = f^*(T, \{ g_{\sigma} \}_{\sigma} )$ in $\aCycle^1(X)_{\QQ}$.
Thus, if we set $y = (1/n)(T, \{ g_{\sigma} \}_{\sigma} )$, then
$x = f^*(y)$ in $\aChow^1(X)_{\QQ}$.

Finally, we assume that there is $y \in \aChow^1(\Spec(O_K))_{\QQ}$
with $x = f^*(y)$ in $\aChow^1(X)_{\QQ}$.
Then, since $\aPic(X)_{\QQ} \otimes \aPic(X)_{\QQ} \to \aChow^2(X)_{\QQ}$
passes through $\aPic(X)_{\QQ} \otimes \aChow^1(X)_{\QQ} \to \aChow^2(X)_{\QQ}$, and
the pairing
\[
\aPic(X)_{\QQ} \otimes \cdots \otimes \aPic(X)_{\QQ} \to \aChow^1(\Spec(O_K))_{\QQ}
\]
is symmetric, we can see
\begin{align*}
\adeg \left( \acherncl_1(\overline{H})^{d-1} \cdot x^2 \right) & =
\adeg \left( \acherncl_1(\overline{H})^{d-1} \cdot (x \cdot f^*(y)) \right) \\
& = \adeg \left( \acherncl_1(\overline{H})^{d-1} \cdot (f^*(y) \cdot x) \right) \\
& = \adeg \left( \acherncl_1(\overline{H})^{d-1} \cdot (f^*(y) \cdot f^*(y)) \right) = 0.
\end{align*}
Thus, we get all assertions of Theorem~\ref{thm:hodge:index:theorem}.
\QED

\begin{Corollary}
\label{cor:hodge:index:type:2}
Let $h \in \aPic_{L^2_1}(X)_{\QQ}$ with $\deg_K(h) > 0$ and
$(h \cdot h)_{\overline{H}} > 0$.
If $x \in \aPic_{L^2_1}(X)_{\QQ}$ and $(h \cdot x)_{\overline{H}} = 0$,
then $(x \cdot x)_{\overline{H}} \leq 0$.
Moreover, if $(x \cdot x)_{\overline{H}} = 0$ and
$f_*(\acherncl_1(H, k)^{d-1} \cdot (h \cdot x)) = 0$, then $x = 0$
in $\aChow^1(X)_{\QQ}$.
\end{Corollary}

\Proof
Let us choose a rational number $t$ with
$\deg_K( x + t h) = 0$.
Then, by Theorem~\ref{thm:hodge:index:theorem},
\[
0 \geq (x + th \cdot x + th)_{\overline{H}} =
(x \cdot x)_{\overline{H}} + t^2 (h \cdot h)_{\overline{H}}.
\]
Thus, $(x \cdot x)_{\overline{H}} \leq 0$.

Next we assume that $(x \cdot x)_{\overline{H}} = 0$ and
$f_*(\acherncl_1(H, k)^{d-1} \cdot (h \cdot x)) = 0$.
Then, in the above inequality, we can see that $t = 0$.
Thus, $\deg_K(x) = 0$ and $(x \cdot x)_{\overline{H}} = 0$.
Therefore, by Theorem~\ref{thm:hodge:index:theorem},
there is  $y \in \aChow^1(\Spec(O_K))_{\QQ}$ with
$x = f^*(y)$ in $\aChow^1(X)_{\QQ}$. Then,
by virtue of Proposition~\ref{prop:projection:formula:chern:ch},
\[
f_* \left(f^*(y) \cdot (\acherncl_1(H, k)^{d} \cdot h) \right) =
y \cdot f_* (\acherncl_1(H, k)^{d} \cdot h)
= \deg_K(h) y.
\]
Thus,
using (1) of Proposition~\ref{prop:assoc:operator},
we have
\begin{align*}
\deg_K(h) y & = f_* \left(f^*(y) \cdot (\acherncl_1(H, k)^{d} \cdot h) \right) =
f_* \left(\acherncl_1(H, k)^{d} \cdot (f^*(y) \cdot h) \right) \\
& = f_* \left( \acherncl_1(H, k)^{d} \cdot (h \cdot f^*(y)) \right) = 0.
\end{align*}
Therefore, $y = 0$. Hence, $x = 0$ in $\aChow^1(X)_{\QQ}$.
\QED

Further, we can give a generalization of \cite[Theorem~A]{MoHI}.

\begin{Corollary}
\label{cor:standard:conj:div}
Let us consider a homomorphism
\[
L : \aChow^p_D(X)_{\QQ} \to \aChow^{p+1}_D(X)_{\QQ}
\]
given by $L(x) = \acherncl_1(\overline{H}) \cdot x$. 
By abuse of notation, the composition of homomorphisms
\[
\aDPic(X)_{\QQ} \longrightarrow \aDChow^1(X)_{\QQ} 
\overset{L^{d-1}}{\longrightarrow} \aDChow^{d}(X)_{\QQ}
\]
is also denoted by $L^{d-1}$.
Then, we have the following.
\begin{enumerate}
\renewcommand{\labelenumi}{(\arabic{enumi})}
\item
$\Ker\left( L^{d-1} : \aDPic(X)_{\QQ} \to \aDChow^{d}(X)_{\QQ}\right)
= \Ker\left(\aPic(X)_{\QQ} \to \aChow^{1}(X)_{\QQ} \right)$.
In particular, if $X$ is regular, then
$L^{d-1} : \aChow^1_D(X)_{\QQ} \to \aChow^{d}_D(X)_{\QQ}$
is injective.

\item
If $x \in \aPic_{L^2_1}(X)_{\QQ}$, $x \not= 0$ in $\aChow^1_{L^2_1}(X)_{\QQ}$, and
$L^d(x) = 0$, then $\adeg(L^{d-1}(x) \cdot x) < 0$.
\end{enumerate}
\end{Corollary}

\Proof
First, let us see (2).
By virtue of (2) of Proposition~\ref{prop:assoc:operator},
\[
L^{d-1}(x) \cdot x = (\acherncl_1(\overline{H})^{d-1} \cdot x) \cdot x
= \acherncl_1(\overline{H})^{d-1} \cdot x^2.
\]
Thus,
applying Corollary~\ref{cor:hodge:index:type:2} to the case 
where $h = \acherncl_1(\overline{H})$,
we have (2). 

\medskip
Next, let us see (1).
It is sufficient to show that if $x \in \aDPic(X)_{\QQ}$ and
$L^{d-1}(x) = 0$, then $x \in \aPic(X)_{\QQ}$ and $x = 0$
in $\aChow^1(X)_{\QQ}$.
Let us choose $x_0 \in \aPic(X)_{\QQ}$ such that
$z(x) = z(x_0)$ and $\omega(x_0)$ is harmonic with respect 
to $c_1(\overline{H})$.
Then, there is a distribution $T$ on $X(\CC)$ with $x = x_0 + a(T)$.
Here, $z(x) \cdot c_1(H)^{d-1} = z(L^{d-1}(x)) = 0$.
Thus, $\omega(x_0) \wedge c_1(\overline{H})^{d-1} = 0$. Therefore,
\begin{align*}
\acherncl_1(\overline{H})^{d-1} \cdot x_0^2 & =
\acherncl_1(\overline{H})^{d-1} \cdot x_0^2 + 
a(\omega(x_0) \wedge c_1(\overline{H})^{d-1} T) \\
& = \acherncl_1(\overline{H})^{d-1} \cdot x_0 \cdot x = x_0 \cdot L^{d-1}(x) = 0.
\end{align*}
Thus, by virtue of Theorem~\ref{thm:hodge:index:theorem},
there is $y \in \aChow^1(\Spec(O_K))_{\QQ}$ with
$x_0 = f^*(y)$ in $\aChow^1(X)_{\QQ}$.
In particular, $\omega(x_0) = 0$.
Therefore, 
\begin{align*}
c_1(\overline{H})^{d-1} \wedge dd^c(T) & =
\omega(\acherncl_1(\overline{H})^{d-1}) \wedge \omega(x_0) +
c_1(\overline{H})^{d-1} \wedge dd^c(T) \\
& = \omega(L^{d-1}(x_0) + a(c_1(\overline{H})^{d-1}T)) =
\omega(L^{d-1}(x)) = 0.
\end{align*}
This implies that $\Delta(T) = 0$, where $\Delta$ is the Laplacian
with respect to $c_1(\overline{H})$.
Hence, using the regularity of solutions of the elliptic operator $\Delta$,
$T$ is represented by the $C^{\infty}$-function $c$
which is constant on each connected component of $X(\CC)$.
In particular, $x \in \aPic^1(X)_{\QQ}$. Thus, by virtue of (2),
we have $x = 0$ in $\aChow^1(X)_{\QQ}$.
\QED

\section{Comparison of intersection numbers via birational morphism}
\label{sec:comp:intersection:numbers:via:birational:morphism}
Let $\mu : Y \to X$ be a birational morphism of
projective arithmetic varieties with $d = \dim X_{\QQ} = \dim Y_{\QQ}$.
We assume that $X$ is normal.
We set
\[
\Pic_{\mu}(Y)_{\QQ} = \{ x \in \Pic(Y)_{\QQ} \mid
\text{$\mu_*(x) \in \Pic(X)_{\QQ}$} \}
\]
and
\[
\aPic_{\mu}(Y)_{\QQ} = \{ x \in \aPic(Y)_{\QQ} \mid z(x) \in 
\Pic_{\mu}(Y)_{\QQ} \},
\]
where $z : \aPic(Y)_{\QQ} \to \Pic(Y)_{\QQ}$ is the homomorphism
forgetting Green functions. Then, the push-forward $\mu_*$ induces
the homomorphism
\[
\mu_* : \aPic_{\mu}(Y)_{\QQ} \to \aBrPic(X)_{\QQ}.
\]
Let us fix an arithmetically ample
Hermitian line bundle $\overline{H} = (H, k)$ on $X$.
Here we define the symmetric bi-linear map
\[
\Delta_{\mu} : \aPic_{\mu}(Y)_{\QQ} \times \aPic_{\mu}(Y)_{\QQ} \to \RR
\]
to be
\[
\Delta_{\mu}(x, y) = \adeg \left( \acherncl_1(\mu^*(\overline{H}))^{d-1} 
\cdot x \cdot y \right) -
\adeg \left( \acherncl_1(\overline{H})^{d-1} \cdot \mu_*(x) \cdot
\mu_*(y) \right).
\]
In this section, we will investigate basic properties of $\Delta_{\mu}$
(cf. Proposition~\ref{prop:Delta:depend:only:delta} and 
\ref{prop:negative:semi:def:comp}) and
give its application (cf. Corollary~\ref{cor:inequality:c2:and:c1:c1}).

\bigskip
Here let us introduce the homomorphism
$\delta_{\mu} : \Pic_{\mu}(Y)_{\QQ} \to \Pic_{\mu}(Y)_{\QQ}$
given by
\[
\delta_{\mu}(x) = x - \mu^*(\mu_*(x)).
\]
By abuse of notation, the composition
\[
\delta_{\mu} \cdot z : \aPic_{\mu}(Y)_{\QQ} \overset{z}{\longrightarrow}
\Pic_{\mu}(Y)_{\QQ} \overset{\delta_{\mu}}{\longrightarrow}
\Pic_{\mu}(Y)_{\QQ}
\]
is also denoted by $\delta_{\mu}$.
First, let us consider the following proposition.

\begin{Proposition}
\label{prop:Delta:depend:only:delta}
If $\delta_{\mu}(x) = \delta_{\mu}(x')$ and $\delta_{\mu}(y) = \delta_{\mu}(y')$
for $x,x',y,y' \in \aPic_{\mu}(Y)_{\QQ}$,
then $\Delta_{\mu}(x, y) = \Delta_{\mu}(x', y')$.
\end{Proposition}

\Proof
First, let us see that
if $z(x) = z(x')$ and $z(y) = z(y')$
for $x,x',y,y' \in \aPic_{\mu}(Y)_{\QQ}$,
then $\Delta_{\mu}(x, y) = \Delta_{\mu}(x', y')$.
For this purpose, it is sufficient to see that
$\Delta_{\mu}(x, a(\phi)) = 0$ for all $C^{\infty}$-functions
$\phi$ on $Y(\CC)$.
First of all,
\[
\adeg \left( \acherncl_1(\mu^*(\overline{H}))^{d-1} 
\cdot x \cdot a(\phi) \right) =
\frac{1}{2} \int_{Y} \phi \omega(x) \wedge \mu^*(\Omega^{d-1}).
\]
On the other hand, we set $x = (D, g)$.
Let $g'$ be a Green function for $\mu_*(D)$.
If we set $\psi = \mu_*(g) - g'$, then $\psi \in \Br(X(\CC))$.
By the definition of the star product and
Remark~\ref{rem:intersection:Br},
\begin{align*}
(\mu_*(D), \mu_*(g)) \cdot (0, \mu_*(\phi)) & =
(0, \omega(g')\mu_*(\phi) + dd^c(\psi) \mu_*(\phi))\\
&  = (0, \omega(\mu_*(g)) \mu_*(\phi)) = (0, \mu_*(\omega(g) \phi)).
\end{align*}
Thus, we have
\[
\adeg \left( \acherncl_1(\overline{H})^{d-1} 
\cdot \mu_*(x) \cdot \mu_*(a(\phi)) \right) =
\frac{1}{2} \int_{Y} \phi \omega(x) \wedge \mu^*(\Omega^{d-1}).
\]
Hence, $\Delta_{\mu}(x, a(\phi)) = 0$.

\medskip
Let us pick up $x, y \in \aPic_{\mu}(Y)_{\QQ}$.
In order to see that $\Delta_{\mu}(x, y)$ depends only on 
$\delta_{\mu}(x)$ and $\delta_{\mu}(y)$,
by replacing $x$ and $y$ by $nx$ and $ny$ for some positive integer $n$,
we may assume that $x, y \in \aPic_{\mu}(Y)$ and
$\mu_*(x), \mu_*(y) \in \aBrPic(X)$.
Let $(L, h_{L})$ and $(Q, h_Q)$ be Hermitian line bundles on $Y$ with
$\acherncl_1(L, h_L) = x$ and $\acherncl_1(Q, h_Q) = y$, and let
$L'$ and $Q'$ be line bundles on $X$ with
$c_1(L') = \mu_*(z(x))$ and $c_1(Q') = \mu_*(z(y))$.
Here we can find Cartier divisors $\Sigma_1$ and $\Sigma_2$
on $Y$ such that $\mu^*(L') \otimes \OO_Y(\Sigma_1) = L$,
$\mu^*(Q') \otimes \OO_Y(\Sigma_2) = Q$, and
$\mu_*(\Sigma_1) = \mu_*(\Sigma_2) = 0$ in $\Div(X)$.
Then, $\delta_{\mu}(x) = \Sigma_1$ and $\delta_{\mu}(y) = \Sigma_2$.
Let $h_{L'}$ and $h_{Q'}$ be $C^{\infty}$-Hermitian metrics of
$L'$ and $Q'$ respectively.
Then, there are $C^{\infty}$-Hermitian metrics
$e_1$ and $e_2$ of $\OO_Y(\Sigma_1)$ and $\OO_Y(\Sigma_2)$
respectively such that
$\mu^*(L', h_{L'}) \otimes (\OO_Y(\Sigma_1), e_1) = (L, h_L)$ and
$\mu^*(Q', h_{Q'}) \otimes (\OO_Y(\Sigma_2), e_2) = (Q, h_Q)$.
Therefore, since 
\[
\begin{cases}
\mu_*(x) = \acherncl_1(L', h_{L'}) +
\mu_*(\acherncl_1(\OO_Y(\Sigma_1), e_1)) \\
\mu_*(y) = \acherncl_1(Q', h_{Q'}) +
\mu_*(\acherncl_1(\OO_Y(\Sigma_2), e_2)),
\end{cases}
\]
using projection formula 
(cf. Proposition~\ref{prop:projection:formula:chern:ch}),
we can easily see that
\[
\Delta_{\mu}(x, y) = \Delta_{\mu}(\acherncl_1(\OO_Y(\Sigma_1), e_1),
\acherncl_1(\OO_Y(\Sigma_2), e_2)).
\]
Thus, combining the first assertion, we have our proposition.
\QED

Before starting the next property, we would like to fix a terminology.
Let $L$ be a line bundle on $Y$, and let $\Gamma$ 
be a prime divisor on $Y$ with $\mu_*(\Gamma) = 0$.
We define $\deg_{\mu}(\rest{L}{\Gamma})$ to be
\[
\deg_{\mu}(\rest{L}{\Gamma}) =
\begin{cases}
\text{the degree of $L$ on the generic fiber of $\Gamma \to \mu(\Gamma)$}
& \text{if $\dim \Gamma - \dim \mu(\Gamma) = 1$} \\
0 & \text{otherwise}
\end{cases}
\]

Let $D_1, \ldots, D_n$ be effective Cartier divisors on $Y$ with
the following properties.
\begin{enumerate}
\renewcommand{\labelenumi}{(\arabic{enumi})}
\item
$\mu_*(D_i) = 0$ for all $i$.

\item
$D_i$ and $D_j$ have no common component for all $i \not= j$.

\item
There are positive integers $a_1, \ldots, a_n$ such that
if we set $D = -\sum_{i=1}^n a_i D_i$, then
$\deg_{\mu}(\rest{\OO_Y(-D)}{\Gamma}) \geq 0$ for all
prime divisors $\Gamma$ in $\Supp(D)$.
\end{enumerate}
Note that if $-D$ is $\mu$-nef, then
(3) is satisfied.

Here we define the subspace $V$ of $\aDiv(X)_{\QQ}$
to be
\[
V = \left\{ (D, g) \in \aDiv(X)_{\QQ} \mid 
\text{$D = \sum_{i=1}^n x_i D_i$ for some $x_1, \ldots, x_n \in \QQ$} 
\right\}.
\]
Then, we have the following proposition.

\begin{Proposition}
\label{prop:negative:semi:def:comp}
$\Delta_{\mu}(x, x') \leq 0$ for any $x, x' \in V$ with $z(x) = z(x')$.
In particular, $\Delta_{\mu}$ is negative semi-definite on $V$.
\end{Proposition} 

\Proof
First of all, note that
\addtocounter{Claim}{1}
\begin{equation}
\label{eqn:prop:negative:semi:def:comp:1}
\adeg \left( \acherncl_1(\overline{H})^{d-1} \cdot \mu_*(D, g) \cdot
\mu_*(D', g') \right) =
\frac{1}{2} \int_Y g' \omega(D, g) \mu^*(\Omega^{d-1})
\end{equation}
because $\mu_*(D, g) = (0, \mu_*(g))$, $\mu_*(D', g') = (0, \mu_*(g'))$
and Remark~\ref{rem:intersection:Br}.
Further, in order to prove our proposition, we may assume
$a_1 = \cdots = a_n = 1$ by replacing $D_i$ by $a_i D_i$.
Let us choose a Green function $g_i$ for each $D_i$.
We set $e_i = (D_i, g_i)$ for $i=1, \ldots, n$, and
$V' = \QQ e_1 + \cdots + \QQ e_n$.
By virtue of Proposition~\ref{prop:Delta:depend:only:delta},
if $\Delta_{\mu}(x, x) \leq 0$ for all $x \in V'$, then
the assertion of our proposition holds.
Here let us consider
the following claim, which is the crucial part of the proof
of our proposition.

\begin{Claim}
\label{claim:prop:negative:semi:def:comp:3}
\begin{enumerate}
\renewcommand{\labelenumi}{(\roman{enumi})}
\item
$\Delta_{\mu}(e_i, e_j) \geq 0$ for all $i \not= j$.

\item
$\Delta_{\mu}(e, e_j) \leq 0$ for all $j$,
where  $e = e_1 + \cdots + e_n$.
\end{enumerate}
\end{Claim}

To prove the above claim, we need the following lemma.

\begin{Lemma}
\label{lem:pos:intersect:num:along:exec:div}
Let $\overline{L}$ be a Hermitian line bundle on $Y$, and
$\Gamma$ a prime divisor on $Y$ with $\mu_*(\Gamma) = 0$.
Let $\nu : \Gamma' \to \Gamma$ be a birational morphism
of projective integral schemes.
We assume that
if $\Gamma$ is horizontal with respect to
$Y \to \Spec(\ZZ)$, then $\nu$ is a generic resolution of
singularities of $\Gamma$; 
otherwise, $\nu = \operatorname{id}_{\Gamma}$.
If $\deg_{\mu}(\rest{L}{\Gamma}) \geq 0$, then
\[
\adeg \left( \acherncl_1(\nu^*(\overline{L})) \cdot 
\acherncl_1(\nu^* \mu^* (\overline{H}))^{d-1} \right) \geq 0.
\]
\rom{(}For the definition of generic resolution of singularities,
see \ref{appendix:projection:formula}.\rom{)}
\end{Lemma}

\Proof
We set $\Sigma = f(\Gamma)$.
Let $\pi : \Sigma' \to \Sigma$ be a proper birational morphism
of projective integral schemes. 
We assume that
if $\Sigma$ is horizontal, then $\pi$ is a generic resolution of
singularities of $\Sigma$; 
otherwise, $\pi = \operatorname{id}_{\Sigma}$.
Changing a model of $\Gamma'$, if necessarily,
we may assume that there is a morphism
$f : \Gamma' \to \Sigma'$ with
$\pi \cdot f  = \mu \cdot \nu$.
\[
\begin{CD}
\Gamma' @>{\nu}>> Y\\
@V{f}VV @VV{\mu}V \\
\Sigma' @>{\pi}>> X
\end{CD}
\]
Thus, using projection formula,
\begin{align*}
\adeg \left( \acherncl_1(\nu^*(\overline{L})) \cdot 
\acherncl_1(\nu^* \mu^* (\overline{H}))^{d-1} \right) & =
\adeg \left( \acherncl_1(\nu^*(\overline{L})) \cdot 
\acherncl_1( f^* \pi^* (\overline{H}))^{d-1} \right) \\
& = \adeg \left( f_* (\acherncl_1(\nu^*(\overline{L}))) \cdot 
\acherncl_1( \pi^* (\overline{H}))^{d-1} \right) \\
& = \deg_{\mu}(\rest{L}{\Gamma})
\adeg \left(  \acherncl_1( \pi^* (\overline{H}))^{d-1} \right) \geq 0.
\end{align*}
\QED

Let us go back to the proof of 
Claim~\ref{claim:prop:negative:semi:def:comp:3}.
Let $D_j = b_1 \Gamma_1 + \cdots + b_s \Gamma_s$
be the irreducible decomposition as cycles.
For each $k$, let $\nu_k : \Gamma'_k \to \Gamma_k$ 
be a proper birational morphism
of projective integral schemes. 
We assume that if $\Gamma_k$ is horizontal with respect to
$Y \to \Spec(\ZZ)$, then $\nu_k$ is a generic resolution of
singularities of $\Gamma_k$.
Then, by Lemma~\ref{lem:formula:restriction:intersection},
\begin{multline*}
\adeg \left( \acherncl_1(\mu^* \overline{H})^{d-1} \cdot
(D_i, g_i) \cdot (D_j, g_j) \right) = \\
\sum_{k=1}^s b_k \adeg \left( 
\acherncl_1(\nu_k^* \mu^* \overline{H})^{d-1} \cdot
\nu_k^* (D_i, g_i) \right) +
\frac{1}{2} \int_{Y(\CC)} g_j \omega(g_i) \mu^{*}(\Omega^{d-1}),
\end{multline*}
where $\Omega = c_1(H, k)$.
Thus, combining the above with \eqref{eqn:prop:negative:semi:def:comp:1},
we can see that
\[
\Delta_{\mu}(e_i, e_j) = \sum_{k=1}^s b_k \adeg \left( 
\acherncl_1(\nu_k^* \mu^* \overline{H})^{d-1} \cdot
\nu_k^* (D_i, g_i) \right).
\]
Here, since $\Gamma_k$ is not a component of $D_i$,
$\deg_{\mu}(\rest{\OO_Y(D_i)}{\Gamma_k}) \geq 0$ for every $k$.
Therefore, by Lemma~\ref{lem:pos:intersect:num:along:exec:div}, 
we get $\Delta_{\mu}(e_i, e_j) \geq 0$.

Let $h$ be a Hermitian metric of $\OO_Y(D)$ with
$e = \acherncl_1(\OO_Y(D), h)$.
In the same way as above,
by using Lemma~\ref{lem:formula:restriction:intersection},
we can see that
\[
\Delta_{\mu}(e, e_j) = \sum_{k=1}^s b_k \adeg \left( 
\acherncl_1(\nu_k^* \mu^* \overline{H})^{d-1} \cdot
\acherncl_1(\nu_k^* (\OO_Y(D), h)) \right).
\]
Here, by our assumption, $\deg_{\mu}(\rest{\OO_Y(-D)}{\Gamma_k}) \geq 0$
for all $k$.
Thus, by Lemma~\ref{lem:pos:intersect:num:along:exec:div},
$\Delta_{\mu}(e, e_j) \leq 0$.

\medskip
Finally, let us see that
$\Delta_{\mu}(x, x) \leq 0$ for all $x \in V'$.
We set $x = x_1 e_1 + \cdots + x_n e_n$.
It is easy to see that
\[
\Delta_{\mu}(x, x) = \sum_i x_i^2 \Delta_{\mu}(e_i, e) - 
\sum_{i < j} (x_i - x_j)^2 \Delta_{\mu}(e_i, e_j).
\]
Thus, Claim~\ref{claim:prop:negative:semi:def:comp:3}
implies $\Delta_{\mu}(x, x) \leq 0$.
\QED

\begin{Corollary}
\label{cor:negative:def:Delta}
We assume that $X$ and $Y$ are $\QQ$-factorial, i.e.,
the natural homomorphisms
$\Div(X)_{\QQ} \to \Cycle^1(X)_{\QQ}$ and
$\Div(Y)_{\QQ} \to \Cycle^1(Y)_{\QQ}$ are isomorphisms.
Then, $\Delta_{\mu}(x, y) \leq 0$ for any $x, y \in \aChow^1(Y)_{\QQ}$
with $\delta_{\mu}(x) = \delta_{\mu}(y)$.
In particular, $\Delta_{\mu}$ is negative semi-definite on
$\aChow^1(Y)_{\QQ}$.
\end{Corollary}

\Proof
Let $\Gamma_1, \ldots, \Gamma_n$ be all prime divisors on $Y$
contracted by $\mu$.
By virtue of Proposition~\ref{prop:Delta:depend:only:delta} and 
\ref{prop:negative:semi:def:comp},
it is sufficient to show that
there are positive rational numbers $a_1, \ldots, a_n$ such that
$-\sum_{i=1}^n a_i \Gamma_i$ is $\mu$-ample.

Let $A$ be an ample line bundle on $Y$.
Then, there are a positive integer $m$ and a section $s$
of $H^0(Y, A^{\otimes m})$ such that
$\Gamma_i$ is not a component of $\zero(s)$ for every $i$.
We set $D = \mu^*(\mu_*(\zero(s)) - \zero(s)$.
(Note that $\mu_*(\zero(s)) \in \Div(X)_{\QQ}$.)
Then, by our choice of $s$, 
$D$ is effective, $-D$ is $\mu$-ample, and $D$ is contracted by $\mu$.
Thus, there are non-negative rational numbers $a_1, \ldots, a_n$
with $D = \sum_{i=1}^n a_i \Gamma_i$.
Here we suppose $a_i = 0$ for some $i$.
Let $F$ be the generic fiber of $\Gamma_i \to \mu(\Gamma_i)$.
Then, $\dim F \geq 1$, $\rest{D}{F}$ is effective, and $-\rest{D}{F}$ is ample.
This is a contradiction. Thus, $a_i > 0$ for all $i$.
\QED

Let $M$ be a $d$-dimensional smooth projective variety over
an algebraically closed field, $E$ a vector bundle bundle of rank $2$
on $M$, and $H$ an ample line bundle on $M$.
It is well known that if $L$ is a rank $1$ saturated subsheaf of $E$,
then
\[
\deg (c_2(E) \cdot c_1(H)^{d-2}) \geq 
\deg (c_1(L) \cdot c_1(E/L) \cdot c_1(H)^{d-2}).
\]
This is very useful formula to estimate the degree of $c_2(E)$.
The following is an arithmetic analogue of the above formula.

\begin{Corollary}
\label{cor:inequality:c2:and:c1:c1}
Let $X$ be a regular projective arithmetic variety
with $d = \dim X_{\QQ}$,
$(E, h)$ a Hermitian vector bundle of rank $2$ on $X$, and
$\overline{H} = (H, k)$ an arithmetically ample Hermitian line
bundle on $X$.
Let $L$ be a saturated subsheaf of $E$ with $\rank L = 1$, and
let $Q = \left( E/L \right)^{**}$.
Let $h_L$ and $h_Q$ be metrics of $L$ and $Q$ induced by
$h$. Then, $h_L$ and $h_Q$ are birationally Hermitian metrics, and
\[
\adeg \left(
\acherncl_1(\overline{H})^{d-1} \cdot \acherncl_2(E, h) \right)
\geq
\adeg \left(
\acherncl_1(\overline{H})^{d-1} \cdot \acherncl_1(L, h_L) \cdot
\acherncl_1(Q, h_Q) \right).
\]
\end{Corollary}

\Proof
First of all, there is an ideal sheaf $I$ on $X$ such that
the image of $E \to Q$ is $Q \otimes I$.
Then, $\codim \Spec(\OO_X/I) \geq 2$.
Let $p_I : X_I \to X$ be the blowing-up by the ideal sheaf $I$, and
$\nu : Y \to X_I$ a generic resolution of singularities of $X_I$.
Then, by our construction,
there is an effective Cartier divisor $\Sigma$ on $Y$
with $I \OO_Y = \OO_Y(-\Sigma)$.
Then, $\mu^*(E) \to \mu^*(Q) \otimes \OO_Y(-\Sigma)$
is surjective and its kernel is
$\mu^*(L) \otimes \OO_Y(\Sigma)$, where $\mu = p_I \cdot \nu : Y \to X$.
Thus, we get an exact sequence
\[
0 \to \mu^*(L) \otimes \OO_Y(\Sigma) \to
\mu^*(E) \to \mu^*(Q) \otimes \OO_Y(-\Sigma) \to 0.
\]
Let $h'_L$ and $h'_Q$ be Hermitian metrics of
$\mu^*(L) \otimes \OO_Y(\Sigma)$ and
$\mu^*(Q) \otimes \OO_Y(-\Sigma)$ induced by $\mu^*(h)$ of $\mu^*(E)$.
Let $\overline{A} = (A, e)$ be an arithmetically ample Hermitian line bundle on $Y$.
Then, by \cite[Proposition~7.3]{MoBGIS}, for all $n > 0$,
\begin{multline*}
\adeg \left(
\acherncl_1(\mu^*(\overline{H}^{\otimes n}) \otimes \overline{A})^{d-1} 
\cdot \acherncl_2(\mu^*(E), \mu^*(h)) \right)
\geq \\
\adeg \left(
\acherncl_1(\mu^*(\overline{H}^{\otimes n}) \otimes \overline{A})^{d-1} 
\cdot \acherncl_1(\mu^*(L) \otimes \OO_Y(\Sigma),\ h'_L) \cdot
\acherncl_1(\mu^*(Q) \otimes \OO_Y(-\Sigma),\  h'_Q) \right).
\end{multline*}
Taking $n \to \infty$ of the above inequality,
we have
\begin{multline*}
\adeg \left(
\acherncl_1(\mu^*(\overline{H}))^{d-1} 
\cdot \acherncl_2(\mu^*(E), \mu^*(h)) \right)
\geq \\
\adeg \left(
\acherncl_1(\mu^*(\overline{H}))^{d-1} 
\cdot \acherncl_1(\mu^*(L) \otimes \OO_Y(\Sigma),\ h'_L) \cdot
\acherncl_1(\mu^*(Q) \otimes \OO_Y(-\Sigma),\ h'_Q) \right).
\end{multline*}
Here, by the projection formula
(cf. Proposition~\ref{prop:projection:formula:chern:ch}),
\[
\adeg \left(
\acherncl_1(\mu^*(\overline{H}))^{d-1} 
\cdot \acherncl_2(\mu^*(E), \mu^*(h)) \right)
=
\adeg \left(
\acherncl_1(\overline{H})^{d-1} 
\cdot \acherncl_2(E, h) \right).
\]
Thus, it is sufficient to show that
\begin{multline*}
\adeg \left(
\acherncl_1(\mu^*(\overline{H}))^{d-1} 
\cdot \acherncl_1(\mu^*(L) \otimes \OO_Y(\Sigma),\ h'_L) \cdot
\acherncl_1(\mu^*(Q) \otimes \OO_Y(-\Sigma),\ h'_Q) \right) \geq \\
\adeg \left(
\acherncl_1(\overline{H})^{d-1} \cdot \acherncl_1(L, h_L) \cdot
\acherncl_1(Q, h_Q) \right).
\end{multline*}
Namely,
\[
\Delta_{\mu}(\acherncl_1(\mu^*(L) \otimes \OO_Y(\Sigma),\ h'_L),\ 
\acherncl_1(\mu^*(Q) \otimes \OO_Y(-\Sigma),\ h'_Q)) \geq 0.
\]
Let $e$ be a Hermitian metric of $\OO_Y(\Sigma)$.
Then, by Proposition~\ref{prop:Delta:depend:only:delta},
\begin{multline*}
\Delta_{\mu}(\acherncl_1(\mu^*(L) \otimes \OO_Y(\Sigma),\ h'_L),\ 
\acherncl_1(\mu^*(Q) \otimes \OO_Y(-\Sigma),\ h'_Q)) = \\
- \Delta_{\mu}(\acherncl_1(\OO_Y(\Sigma),e), \acherncl_1(\OO_Y(\Sigma),e)).
\end{multline*}
Therefore, by Proposition~\ref{prop:negative:semi:def:comp},
it suffices to show that $-\Sigma$ is $\mu$-nef.
This is obvious because
$I\OO_{X_I}$ is $\mu_I$-ample and $\OO_Y(-\Sigma) = \nu^*(I\OO_{X_I})$.
\QED

\section{Bogomolov's instability for rank 2 bundles}

Let $K$ be a number field, and $O_K$ the ring of integers of $K$.
Let $f : X \to \Spec(O_K)$ be a regular projective arithmetic variety
with the geometrically irreducible generic fiber and $d = \dim X_{K}$.
Let $\overline{H} = (H, k)$ be an arithmetically ample Hermitian
line bundle on $X$. Let
\[
(\ \cdot \ )_{\overline{H}} : \aBrChow^1(X)_{\QQ} \times \aBrChow^1(X)_{\QQ}
\to \RR
\quad\text{and}\quad
\deg_K : \aBrChow^1(X)_{\QQ} \to \QQ
\]
be homomorphisms given in \S\ref{sec:hodge:index:thm}.
Here, we set
\[
\poscone(X; \overline{H}) = \{
x \in \aBrChow^1(X)_{\QQ} \mid
\text{$(x \cdot x)_{\overline{H}} > 0$ and $\deg_K(x) > 0$} \}
\]
and
\[
\wposcone(X; \overline{H}) = \{
x \in \aBrChow^1(X)_{\QQ} \mid
\text{$(x \cdot y)_{\overline{H}} > 0$ for all $y \in \poscone(X; \overline{H})$}
\}.
\]
By virtue of the Hodge index theorem 
(cf. Theorem~\ref{thm:hodge:index:theorem}), 
we have the following in the same way as \cite[\S1]{MoBGUN}.

\begin{Proposition}
\label{prop:poscone:in:wposcone}
$\poscone(X; \overline{H}) \subset \wposcone(X; \overline{H})$.
\end{Proposition}

Let $\overline{E} = (E, h)$ be a Hermitian vector bundle of rank $2$ on $X$. 
Let $L$ be a saturated subsheaf of $E$ with $\rank L = 1$.
Since $X$ is regular and $L$ is reflexive, $L$ is an invertible
sheaf. Let $h_L$ be the metric induced by $h$ of $E$.
Then, $h_L$ is a birationally Hermitian metric. In this notation,
we have the following theorem.

\begin{Theorem}
\label{them:B:instability:rank:2}
If $\adeg\left( \acherncl_1(\overline{H})^{d-1} \cdot 
\left( 4 \acherncl_2(\overline{E}) - \acherncl_1(\overline{E})^2 \right)
\right) < 0$, there is a saturated rank $1$ subsheaf $L$ of $E$
with
\[
2 \acherncl_1(L, h_L) - \acherncl_1(E, h) \in
\poscone(X; \overline{H}).
\]
\end{Theorem}

\Proof
First of all, by virtue of \cite{MoBGI},
$E_{\overline{K}}$ is not $\mu$-semistable with respect to
$H_{\overline{K}}$. Thus, there is the destabilizing subsheaf $L'$
of $E_{\overline{K}}$. Using the uniqueness of $L'$,
in the same way as \cite[Claim~3.2]{MoBGUN},
we can see that $L'$ is defined over $K$.
Hence, we can find a saturated subsheaf $L$ of $E$
with $L_K = L'$.
Let $Q = (E/L)^{**}$, and let $h_L$ and $h_Q$ be
birationally Hermitian metrics of $L$ and $Q$ induced by $h$.
Then, by Corollary~\ref{cor:inequality:c2:and:c1:c1},
\[
\adeg \left(
\acherncl_1(\overline{H})^{d-1} \cdot \acherncl_2(E, h) \right)
\geq
\adeg \left(
\acherncl_1(\overline{H})^{d-1} \cdot \acherncl_1(L, h_L) \cdot
\acherncl_1(Q, h_Q) \right).
\]
Therefore, since $\acherncl_1(E, h) = \acherncl_1(L, h_L) + \acherncl_1(Q, h_Q)$,
\begin{align*}
\left( (2\acherncl_1(L, h_L) - \acherncl_1(E, h))^2 \right)_{\overline{H}}
& = (\acherncl_1(E, h)^2)_{\overline{H}} -
4 (\acherncl_1(L, h_L) \cdot \acherncl_1(Q, h_Q))_{\overline{H}} \\
& \geq
\adeg\left( \acherncl_1(\overline{H})^{d-1} \cdot 
\left( \acherncl_1(\overline{E})^2 - 4 \acherncl_2(\overline{E}) \right)
\right) > 0.
\end{align*}
Thus, we get our theorem.
\QED

Let us fix $\lambda \in \wposcone(X; \overline{H})$.
Let $\overline{E} = (E, h)$ be a Hermitian vector bundle on $X$ of
rank $2$.
$\overline{E}$ is said to be {\em arithmetically semistable with
respect to $\lambda$} if, for any saturated rank $1$ subsheaves $L$ of $E$,
\[
(\acherncl_1(L, h_L) \cdot \lambda)_{\overline{H}} \leq
\frac{(\acherncl_1(E, h) \cdot \lambda)_{\overline{H}}}{2}.
\]
With notation as above, we have the following corollary.

\begin{Corollary}
\label{cor:BG:ineq:rank:2:bundle}
If $\overline{E}$ is arithmetically semistable with respect to
$\lambda$, then,
\[
\adeg\left( \acherncl_1(\overline{H})^{d-1} \cdot 
\left( 4 \acherncl_2(\overline{E}) - \acherncl_1(\overline{E})^2 \right)
\right) \geq 0.
\]
\end{Corollary}

\Proof
If $\adeg\left( \acherncl_1(\overline{H})^{d-1} \cdot 
\left( 4 \acherncl_2(\overline{E}) - \acherncl_1(\overline{E})^2 \right)
\right) < 0$, then,
by Theorem~\ref{them:B:instability:rank:2},
there is a saturated subsheaf $L$ of $E$ with
\[
2 \acherncl_1(L, h_L) - \acherncl_1(E, h) \in
\poscone(X; \overline{H}).
\]
Thus,
\[
(2 \acherncl_1(L, h_L) - \acherncl_1(E, h) \cdot \lambda)_{\overline{H}} > 0.
\]
which contradicts to the semistability of $\overline{E}$.
\QED

\renewcommand{\thesection}{Appendix \Alph{section}}
\renewcommand{\theTheorem}{\Alph{section}.\arabic{Theorem}}
\renewcommand{\theClaim}{\Alph{section}.\arabic{Theorem}.\arabic{Claim}}
\renewcommand{\theequation}{\Alph{section}.\arabic{Theorem}.\arabic{Claim}}
\setcounter{section}{0}

\section{Projection formula for the arithmetic Chern Character}
\label{appendix:projection:formula}
In this section, we will show the following projection formula
in a general context.

\begin{Proposition}
\label{prop:projection:formula:chern:ch}
Let $X \to Y$ be a proper morphism of arithmetic varieties,
$\overline{E} = (E, h)$ a Hermitian vector bundle on $Y$, and
$z$ an arithmetic $D$-cycle on $X$. Then,
\[
f_* (\achernch(f^* \overline{E}) \cdot z) = 
\achernch(\overline{E}) \cdot f_* z.
\]
\end{Proposition}

\Proof
For the proof of the projection formula above, 
we need the following two lemmas.
The proof of these lemmas can be found in 
\cite[Proposition~2.4.1 and Proposition~2.4.2]{KMRB}.
Here we fix notation.
Let $Z$ be a quasi-projective integral scheme over $\ZZ$.
Then, by virtue of Hironaka's resolution of singularities \cite{Hiro},
there is a proper birational morphism $\mu : Z' \to Z$
of quasi-projective integral schemes over $\ZZ$ such that
$Z'_{\QQ}$ is non-singular.
The above $\mu : Z' \to Z$ is called a {\em generic resolution of
singularities of $Z$}.

\begin{Lemma}
\label{lem:projection:formula:line:bundle}
Let $f : X \to Y$ be a proper morphism of arithmetic varieties. 
Let $(L, h)$ be a Hermitian line bundle on $Y$, and
$z \in \aDChow^p(X)$. Then
\[
f_*(\acherncl_1 (f^*L, f^*h) \cdot z) = \acherncl_1(L, h) \cdot f_*(z).
\]
\end{Lemma}

\begin{Lemma}
\label{lem:formula:restriction:intersection}
Let $X$ be a arithmetic variety, and
$\overline{L}_1 = (L_1, h_1), \ldots, \overline{L}_{n} = (L_{n}, h_{n})$
be Hermitian line bundles on $X$.
Let $(Z, g)$ be an arithmetic $D$-cycle on $X$, and
$Z = a_1 Z_1 + \cdots + a_r Z_r$ the irreducible decomposition as cycles.
For each $i$, let $\tau_i : Z'_i \to Z_i$ be a proper birational morphism
of quasi-projective integral schemes. We assume that if $Z_i$ is horizontal with respect to
$X \to \Spec(\ZZ)$, then $\tau_i$ is a generic resolution of
singularities of $Z_i$. Then, we have
\[
\acherncl_1(\overline{L}_1) \cdots \acherncl_1(\overline{L}_n) \cdot (Z, g) =
\sum_{i=1}^r a_i {\mu_i}_* \left( 
\acherncl_1(\mu_i^* \overline{L}_1) \cdots \acherncl_1(\mu_i^* \overline{L}_n) 
\right) +
a(c_1(\overline{L}_1) \wedge \cdots \wedge c_1(\overline{L}_n) \wedge g)
\]
in $\aDChow^*(X)_{\QQ}$,
where $\mu_i$ is the composition of
$Z'_i \overset{\tau_i}{\longrightarrow} Z_i \hookrightarrow X$
for each $i$.
\end{Lemma}

Let us start the proof of Proposition~\ref{prop:projection:formula:chern:ch}.
We will prove this proposition by induction on $r = \rank E$.
If $r = 1$, then
\[
\achernch(f^* \overline{E}) = 
\sum_{n \geq 0} \frac{1}{n!} \acherncl_1(f^* \overline{E})^n
\quad\text{and}\quad
\achernch(\overline{E}) = 
\sum_{n \geq 0} \frac{1}{n!} \acherncl_1(\overline{E})^n .
\]
Thus, our proposition is a consequence of
Lemma~\ref{lem:projection:formula:line:bundle}.
Thus, we may assume that $r > 1$.
Moreover, if $z = (0, T)$, then
\begin{align*}
f_* (\achernch(f^* \overline{E}) \cdot z) & =
f_* (0, \operatorname{ch}(f^* \overline{E}) \wedge T)
= (0, f_*( f^* \operatorname{ch}(\overline{E}) \wedge T)) \\
& = (0, \operatorname{ch}(\overline{E}) \wedge f_* T) =
\achernch(\overline{E}) \cdot f_* z.
\end{align*}
Hence, we may further assume that $z$ is an usual arithmetic cycle,
i.e., $z \in \aChow^p(X)$ for some $p$.

Let $\pi : P = \Proj(\oplus_{n \geq 0} \Sym^n(E)) \to Y$
and $\nu : Q = \Proj(\oplus_{n \geq 0} \Sym^n(f^* E)) \to X$
be the projective bundles of $E$ and $f^* E$, and
let $\OO_P(1)$ and $\OO_Q(1)$ be the tautological line
bundles of $P$ and $Q$ respectively.
We set the induced morphisms as following diagram.
\[
\begin{CD}
X @<{\nu}<< Q \\
@V{f}VV @VV{g}V \\
Y @<<{\pi}< P
\end{CD}
\]
We give $\OO_Q(1)$ the Hermitian metric
induced from $\nu^* f^* \overline{E}$.
Since $\nu : Q \to X$ is smooth,
we can consider the pull-back $\nu^*(z)$ of $z$.
Here we claim the following.

\addtocounter{Claim}{1}
\begin{Claim}
\label{claim:puch:nu:z:OO}
$\nu_* (\acherncl_1(\overline{\OO_Q(1)})^{r-1} \cdot \nu^*(z)) = z$.
\end{Claim}

\Proof
Let $(Z, g)$ be a representative of $z$.
Clearly, we may assume that $Z$ is integral.
If $Z$ is vertical, then our assertion is trivial.
So we may assume that $Z$ is horizontal.
Let $\mu : Z' \to Z$ be a generic resolution of singularities
of $Z$. Let $\nu' : T = 
\Proj(\oplus_{n \geq 0} \Sym^n(\mu^* f^*(E))) \to Z'$
be the projective bundle of $\mu^* f^* (E)$, 
and $\OO_T(1)$ the tautological line bundle
on $T$. We set the induced morphism as follows.
\[
\begin{CD}
T @>{\mu'}>> Q \\
@V{\nu'}VV @VV{\nu}V \\
Z' @>>{\mu}> X
\end{CD}
\]
Then, $\mu'$ gives rise to a generic resolution of singularities
of $\nu^*(Z)$. Thus, by virtue of
Lemma~\ref{lem:formula:restriction:intersection},
\[
\acherncl_1(\overline{\OO_Q(1)})^{r-1} \cdot \nu^*(z) =
\mu'_* (\acherncl_1(\overline{\OO_T(1)})^{r-1}) +
a(c_1(\overline{\OO_Q(1)})^{r-1} \wedge \nu^*(g)).
\]
Here since
\[
\nu_* \mu'_* (\acherncl_1(\overline{\OO_T(1)})^{r-1})
= \mu_* \nu'_* (\acherncl_1(\overline{\OO_T(1)})^{r-1})
= \mu_*([Z']) = (Z, 0)
\]
and $\nu_* (c_1(\overline{\OO_Q(1)})^{r-1} \wedge \nu^*(g)) = g$,
we can see that
\[
\nu_* (\acherncl_1(\overline{\OO_Q(1)})^{r-1} \cdot \nu^*(z)) =
(Z, 0) + a(g) = (Z, g).
\]
Hence, we get our claim.
\QED

Let us go back to the proof of 
Proposition~\ref{prop:projection:formula:chern:ch}.
We set $\beta = \acherncl_1(\overline{\OO_Q(1)})^{r-1} \cdot \nu^*(z)$.
Then, by the above claim, $\nu_*(\beta) = z$.
Thus, since $\nu$ is smooth,
using \cite[6 of Theorem~4]{GSRR},
\[
f_* \nu_* (\achernch(\nu^* f^* \overline{E}) \cdot \beta)
= f_* (\achernch(f^* \overline{E}) \cdot z).
\]
On the other hand,
let $F_P$ be the kernel of the natural homomorphism
$\pi^*(E) \to \OO_P(1)$.
We give $F_P$ and $\OO_P(1)$ the metrics induced
from the metric of $\pi^* \overline{E}$,
so that we get an exact sequence of Hermitian vector bundles:
\[
\mathcal{E}_Q : 0 \to \overline{F_P} \to \pi^* \overline{E} \to
\overline{\OO_P(1)} \to 0.
\]
Then,
\[
\achernch(\pi^* \overline{E}) =
\achernch(\overline{F_P}) + \achernch(\overline{\OO_P(1)}) -
a(\widetilde{\operatorname{ch}}(\mathcal{E}_P)),
\]
which implies
\[
\achernch(\nu^* f^* \overline{E}) = \achernch(g^* \pi^* \overline{E})
= \achernch(g^* \overline{F_P}) + \achernch(g^* \overline{\OO_P(1)}) -
a(g^* \widetilde{\operatorname{ch}}(\mathcal{E}_P)).
\]
Thus, since the rank of $F_P$ is less than $r$ and
$\pi$ is smooth,
using hypothesis of induction and \cite[6 of Theorem~4]{GSRR},
we have
\begin{align*}
f_* \nu_* (\achernch(\nu^* f^* \overline{E}) \cdot \beta) & =
\pi_* g_* \left(
(\achernch(g^* \overline{F_P}) + \achernch(g^* \overline{\OO_P(1)}) -
a(g^* \widetilde{\operatorname{ch}}(\mathcal{E}_P))) \cdot \beta
\right) \\
& =
\pi_* \left(
(\achernch(\overline{F_P}) + \achernch(\overline{\OO_P(1)}) -
a(\widetilde{\operatorname{ch}}(\mathcal{E}_P))) \cdot g_* \beta
\right) \\
& = \pi_* ( \achernch(\pi^* \overline{E}) \cdot g_* \beta ) 
= \achernch(\overline{E}) \cdot \pi_* g_* \beta.
\end{align*}
Therefore,
\[
f_* (\achernch(f^* \overline{E}) \cdot z) =
f_* \nu_* (\achernch(\nu^* f^* \overline{E}) \cdot \beta) =
\achernch(\overline{E}) \cdot f_* z
\]
because $\pi_* g_* \beta = f_* \nu_* \beta = f_* z$.
\QED

\section{Zariski's lemma for integral scheme}

Let $R$ be a discrete valuation ring, and
$f : Y \to \Spec(R)$ a flat and projective integral scheme over $R$.
Let $\eta$ be the generic point of $\Spec(R)$ and
$o$ the special point of $\Spec(R)$.
We assume that the genetic fiber $Y_{\eta}$ of $f$ is geometrically
reduced and irreducible. Let $Y_o$ be the special fiber of $f$, i.e.,
$Y_o = f^*(o)$.
Let us consider a paring
\[
\Pic(Y) \otimes \Chow^p(Y_o) \to \Chow^{p+1}(Y_o)
\]
given by the composition of homomorphisms
\[
\Pic(Y) \otimes \Chow^p(Y_o) \to \Pic(Y_o) \otimes \Chow^p(Y_o) \to
\Chow^{p+1}(Y_o).
\]
We denote by $x \cdot z$ the image of $x \otimes z$ by the above 
homomorphism.
Let $D$ be a Cartier divisor on $Y$, and
$Z$ a cycle of codimension $p$ on $Y_o$, i.e., $Z \in Z^p(Y_o)$.
We assume that $D$ and $Z$ intersect properly.
Let $s$ be a rational section of $\OO_Y(D)$ with $\zero(s) = D$, and
let $Z = a_1 Z_1 + \cdots + a_r Z_r$ be the irreducible decomposition
as cycles. Then, since $\rest{s}{Z_i}$ is a rational section of $\OO_{Z_i}(D)$,
we define $D \sqcap Z \in Z^{p+1}(Y_o)$ to be
\[
D \sqcap Z = a_1 \zero(\rest{s}{Z_i}) + \cdots + a_r \zero(\rest{s}{Z_r}).
\]
Then, the class of $D \sqcap Z$ is equal to 
$\OO_Y(D) \cdot (\text{the class of $Z$})$.
Moreover, for a Cartier divisor $D$ on $Y$,
the associated cycle of $D$ is denoted by $[D]$, which is an element
of $Z^1(Y)$.
Let us consider the following subgroup $F_c(Y)$ of $Z^0(Y_o)$:
\[
F_c(Y) = \{ x \in Z^0(Y_o) \mid
\text{$x = [D]$ for some Cartier divisor $D$ on $Y$} \}.
\]
For a Cartier divisor $D$ on $Y$ with $[D] \in F_c(Y)$, and
$y \in F_c(Y)$, $D \cdot y$ depend only on $[D]$.
For, if $D'$ is a Cartier divisor on $Y$ with $[D'] = [D]$, and
$E$ is a Cartier divisor on $Y$ with $y = [E]$, then,
by \cite[Theorem~2.4]{Fu},
\[
D \cdot y = E \cdot [D] = E \cdot [D'] = D' \cdot y.
\]
Thus, we can define a bi-linear map
\[
q : F_c(Y) \times F_c(Y) \to \Chow^1(Y_o)
\]
by $q([D], y) = D \cdot y$.
Moreover, \cite[Theorem~2.4]{Fu} says us that $q$ is symmetric, i.e.,
$q(x, y) = q(y, x)$ for all $x, y \in F_c(Y)$.
Let $H$ be an ample line bundle on $Y$. Using $q$ and $H$, we have
a quadratic form $Q_H$ on $F_c(Y)$ given by
\[
Q_H(x, y) = \deg ( H^{d-1} \cdot q(x, y) ),
\]
where $d = \dim Y_{\eta}$.
Then, we have the following Zariski's lemma on integral schemes,
which is well known if $X$ is regular.

\begin{Lemma}[Zariski's lemma for integral scheme]
\label{lem:int:on:fibers}
\begin{enumerate}
\renewcommand{\labelenumi}{(\arabic{enumi})}
\item
$Q_H([Y_o], x) = 0$ for all $x \in F_c(Y)_{\QQ}$.

\item
$Q_H(x, x) \leq 0$ for any $x \in F_c(Y)_{\QQ}$.

\item
$Q_H(x, x) = 0$ if and only if
$x \in \QQ \cdot [Y_o]$.
\end{enumerate}
\end{Lemma}

\Proof
(1):\quad
This is obvious because $\OO_Y(Y_o) \simeq \OO_Y$.

(2) and (3):\quad
If $x \in \QQ \cdot [Y_o]$, then by (1), $Q_H(x, x) = 0$.
Thus, it is sufficient to prove that
(a) $Q_H(x, x) \leq 0$ for any $x \in F_c(Y)_{\QQ}$, and that
(b) if $Q_H(x, x) = 0$, then $x \in \QQ \cdot [Y_o]$.
For this purpose, we may assume that
$x \in F_c(Y)$, i.e., $x = [D]$ for some Cartier divisor $D$ on $Y$.
We prove (a) and (b) by induction on $d$.
If $d=1$ and $Y$ is regular, the lemma follows from the following sublemma.

\begin{Sublemma}
Let $V$ be a finite dimensional vector space over $\RR$, and
$Q$ a quadratic form on $V$. We assume that there are $e \in V$ and
a basis $\{ e_1, \ldots, e_n \}$ with the following properties:
\begin{enumerate}
\renewcommand{\labelenumi}{(\roman{enumi})}
\item
If we set $e = a_1 e_1 + \cdots + a_n e_n$, then
$a_i > 0$ for all $i$.

\item
$Q(x, e) \leq 0$ for all $x \in V$.

\item
$Q(e_i, e_j) \geq 0$ for all $i \not= j$.

\item
If we set $S = \{ (i, j) \mid \text{$i \not= j$ and $Q(e_i, e_j) > 0$} \}$,
then, for any $i \not= j$, there is a sequence $i_1, \ldots, i_l$ such that
$i_1 = i$, $i_l = j$, and $(i_t, i_{t+1}) \in S$ for all $1 \leq t < l$.
\end{enumerate}
Then, $Q(x, x) \leq 0$ for all $x \in V$.
Moreover, if $Q(x, x) = 0$ for some $x \not= 0$, then
$x \in \RR e$ and $Q(y, e) = 0$ for all $y \in V$.
\end{Sublemma}

\Proof
Replacing $e_i$ by $a_i e_i$, we may assume that $a_1 = \cdots = a_n = 1$.
If we set $x = x_1 e_1 + \cdots + x_n e_n$,
then, by an easy calculation, we can show
\[
Q(x, x) = \sum_i x_i^2 Q(e_i, e) - \sum_{i < j} (x_i - x_j)^2 Q(e_i, e_j).
\]
Thus, we can easily see our assertions.
\QED

Let us go back to the proof of Lemma~\ref{lem:int:on:fibers}.
Let $\mu : Y' \to Y$ be a proper birational morphism
of projective integral schemes over $R$. We assume that
if $d \geq 2$, then $\mu$ is finite.
Here we claim that if (a) and (b) hold for $Y'$, then they hold for $Y$.
Note that if $d=1$, then the lemma does not involve $H$.
By virtue of projection formula (cf. \cite[(c) of Proposition~2.4]{Fu}),
\[
\deg ( \mu^*(H)^{d-1} \cdot \OO_{Y'}(\mu^*(D)) \cdot [\mu^*(D)] ) =
\deg ( H^{d-1} \cdot \OO_Y(D) \cdot [D] ).
\]
Thus, if $Q_H([\mu^*(D)], [\mu^*(H)]) \leq 0$,
then $Q_H([D], [D]) \leq 0$.
Moreover, if there is a rational number $\alpha$ such that
$[\mu^*(D)] = \alpha [Y'_o]$, then
$[\mu^*(D)] = \alpha [\mu^*(Y_o)]$. Thus, taking the push-forward $\mu_*$,
we can see that $[D]= \alpha [Y_o]$ in $\Cycle^1(Y)_{\QQ}$.
Hence, we get our claim.

By the above claim, considering the normalization of $Y$,
we may assume that $Y$ is normal.
Moreover, if $d=1$, 
there is a resolution of singularities $\mu : Y' \to Y$ of $Y$
(cf. \cite{Lip}).
Thus it holds for $d = 1$.
Hence we may assume $d \geq 2$.

Let $(Y_o)_{red} = \Gamma_1 \cup \cdots \cup \Gamma_l$ be the irreducible decomposition
of $(Y_o)_{red}$ and $I$ the defining ideal of $(Y_o)_{red}$.
Since $H$ is ample, there is a positive integer $m$ such that
$H^{\otimes m}$ is very ample and $H^1(Y, H^{\otimes m} \otimes I) = 0$.
Thus, $H^0(Y, H^{\otimes m}) \to H^0((Y_o)_{red}, \rest{H^{\otimes m}}{(Y_o)_{red}})$
is surjective. 
Here note that $\codim(Y \setminus \Sing(Y)) \geq 2$ because $Y$ is normal.
Hence, there is a section $s_0$ of $H^0(X, H^{\otimes m})$ such that
$\rest{s_0}{\Gamma_i} \not= 0$
for every $i$, and that $\zero(\rest{s_0}{(Y_o)_{red}})$ intersects with
$(\Gamma_i \cap \Sing(Y))_{red}$ and
$(\Gamma_i \cap \Gamma_j)_{red}$ properly
for all $i \not= j$.
Let $t$ be a element of $R$ such that $t$ is a generator of the maximal ideal of $R$.
Since $H^0(Y, H^{\otimes m})$ is a free $R$-module, there is a basis
$\{ e_1, \ldots, e_n \}$ of $H^0(Y, H^{\otimes m})$ as $R$-module.
Then, there are $c_1, \ldots, c_n \in R$ with $s_0 = c_1 e_1 + \cdots + c_n e_n$.
For each $a_1, \ldots, a_r$ of $R$, 
let us consider the following element $s$ of 
$H^0(Y, H^{\otimes m})$;
\[
s = s_0 + t (a_1 e_1 + \cdots + a_n e_n) = (c_1 + t a_1) e_1 + \cdots + (c_n + t a_n) e_n.
\]
Since $\#(R) = \infty$, it is easy to see that the set
\[
\{ (c_1 + t a_1, \ldots, c_n + t a_n) \mid a_1, \ldots, a_n \in R \}
\]
is Zariski dense in $\mathbb{A}^n(K)$, where $K$ is the quotient field of $R$.
Thus, we can find $a_1, \ldots, a_n \in R$ such that
$\zero(s_{\eta})$ is geometrically reduced and 
irreducible divisor on $X_{\eta}$
(cf. \cite[Theorem~6.10]{JB}).
Let $[\zero(s)] = S + T$ be the decomposition as cycles such that $S$ is horizontal and
$T$ is vertical with respect to $f$.
Then, by our choice of $s$, 
$\rest{s}{\Gamma_i} = \rest{s_0}{\Gamma_i}$ for all $i$. Thus,
$T = 0$ and $S$ is integral.
Hence, by the proof of \cite[Theorem~2.4]{Fu}, we can see that $\zero(s) \sqcap [D] = D \sqcap S$
in $Z^2(Y)$. Therefore, 
if we set $H_S = \rest{H}{S}$ and $D_S = \rest{D}{S}$, then
\begin{align*}
\deg ( H_S^{d-2} \cdot \OO_S(D_S) \cdot [D_S] ) & =
\deg ( H^{d-2} \cdot \OO_Y(D) \cdot (D \sqcap S) ) \\
& = \deg (  H^{d-2} \cdot \OO_Y(D) \cdot (\zero(s) \sqcap [D]) ) \\
& = \deg (  H^{d-2} \cdot \OO_Y(D) \cdot H^{\otimes m} \cdot [D]) ) \\
& = m \deg (  H^{d-1} \cdot \OO_Y(D) \cdot  [D]) ).
\end{align*}
Thus, by hypothesis of induction, $Q_H([D], [D]) \leq 0$.

Further, we assume that $Q_H([D], [D]) = 0$.
Then, $\deg ( H_S^{d-2} \cdot \OO_S(D_S) \cdot [D_S] ) = 0$.
Thus, by hypothesis of induction,
there is a rational number $\alpha$ with
$D \sqcap S = \alpha (Y_o \sqcap S)$ in $Z^2(Y)_{\QQ}$.
We set $[D] = \sum_i \alpha_i \Gamma_i$ and $[Y_o] = \sum_i \beta_i \Gamma_i$
as cycles.
Moreover, we set
\[
Y^0 = Y
\setminus \left( \Sing(Y) \cup \bigcup_{i \not= j} (\Gamma_i \cap \Gamma_j) \right),
\]
$D^0 = D \cap Y^0$, $Y_o^0 = Y_o \cap Y^0$,
$S^0 = S \cap Y^0$, and $\Gamma_i^0 = \Gamma_i \cap Y^0$ for each $i$.
Then,
\[
D^0 \sqcap S^0 = \sum_i \alpha_i (\Gamma_i^0 \sqcap S^0)
\quad\text{and}\quad
Y_o^0 \sqcap S^0 = \sum_i \beta_i (\Gamma_i^0 \sqcap S^0)
\]
in $Z^2(Y^0)$.
Thus,
\[
\sum_i \alpha_i (\Gamma_i^0 \sqcap S^0) = \sum_i \alpha \beta_i (\Gamma_i^0 \sqcap S^0)
\]
in $Z^2(Y^0)_{\QQ}$.
Here $H^{\otimes m}$ is very ample and
$\zero(\rest{s_0}{(Y_o)_{red}})$ intersects with
$(\Gamma_i \cap \Sing(Y))_{red}$ and $(\Gamma_i \cap \Gamma_j)_{red}$ properly
for all $i \not= j$.
Therefore,
$\Gamma_i^0 \sqcap S^0 \not= 0$ for all $i$, and
$\Gamma_i^0 \sqcap S^0$ and $\Gamma_j^0 \sqcap S^0$ have no common component
for all $i \not= j$.
Thus, we have $\alpha_i = \alpha \beta_i$ for all $i$.
Hence $[D] = \alpha [Y_o]$ in $Z^1(Y)_{\QQ}$.
\QED

\end{document}